\documentclass[10pt,reqno]{amsart}
   \usepackage[centertags]{amsmath}
   \usepackage{amsfonts}
   \usepackage{amsmath}
   \usepackage{amssymb,mathtools}
   \usepackage{amsthm}
   \usepackage{newlfont}
   \usepackage[latin1]{inputenc}%%%%para que reconozca las tildes.

\usepackage{caption}
\usepackage{subcaption}

\usepackage{multirow, bigdelim}

\usepackage{epsfig}
\usepackage{pst-all}
\usepackage{pstricks}

\usepackage{pst-pdgr}

\usepackage{caption}

\usepackage{color}

\usepackage{graphicx}

\usepackage{bm}

\newcommand\D{\displaystyle}

\allowdisplaybreaks[3]
\theoremstyle{plain}
\newtheorem{theorem}{Theorem}[section]

\theoremstyle{definition}

\newtheorem{example}[theorem]{Example}

\theoremstyle{remark}
\newtheorem{remark}[theorem]{Remark}
\newtheorem*{remark*}{Remark}

\numberwithin{equation}{section}

\title[Examples of spectral analysis of bilateral BDP]{Spectral analysis of bilateral birth-death \\ processes: some new explicit examples}

\author{Manuel D. de la Iglesia}
\address{Manuel D. de la Iglesia\\
Instituto de Matem\'aticas \\
Universidad Nacional Aut\'onoma de M\'exico \\
Circuito Exterior, C.U.\\
04510, Ciudad de M\'exico, Mexico.}
\email{mdi29@im.unam.mx}
\thanks{This work was partially supported by PAPIIT-DGAPA-UNAM grant IN104219 (M\'exico) and CONACYT grant A1-S-16202 (M\'exico).}

\date{\today}



\begin{document}

\maketitle

\begin{abstract}

We consider the spectral analysis of several examples of bilateral birth-death processes and compute explicitly the spectral matrix and the corresponding orthogonal polynomials. We also use the spectral representation to study some probabilistic properties of the processes, like recurrence, the invariant distribution (if it exists) or the probability current.

\end{abstract}

\section{Introduction}

Birth-death processes belong to an important sub-class of continuous-time Markov chains. They are characterized by the property that the only possible transitions are between neighboring states. Birth-death processes are frequently found in many areas of science like biology, genetics, ecology, physics, mathematical finance, queuing and communication systems, epidemiology or chemical reactions (see \cite{PL} for an extensive survey about the subject). These processes are usually defined on the state space of nonnegative integers $\mathbb{N}_0$. However, there are some situations in physics, chemistry or engineering where the state space is the set of all integers $\mathbb{Z}$ (see \cite{Con,dCIM,HP,PL,TB}). These processes are usually known as \emph{bilateral birth-death processes}, although some other names like unrestricted birth-death chains or double-ended systems are also found in the literature. They were firstly studied by W.E. Pruitt in \cite{Pru} (see also \cite{PruT}), following the pioneering works of S. Karlin and J. McGregor for birth-death processes on $\mathbb{N}_0$ (see \cite{KMc2,KMc3,KMc6}). The main tool used in the previous papers is the \emph{spectral theorem} applied to the infinitesimal operator associated with a birth-death process, which is a tridiagonal or Jacobi matrix. The spectral analysis of this kind of operators is related with the theory of orthogonal polynomials and it provides an integral representation of the transition probability functions $P_{ij}(t)$, usually called the \emph{Karlin-McGregor formula} (see \eqref{3FKMcGN} below). In the case of bilateral birth-death processes we need to apply the spectral theorem \emph{three times} and we can accommodate these measures in a $2\times2$ matrix which is usually called the \emph{spectral matrix} (see \eqref{3FKMcGZ} and \eqref{3FKMcGZ2} below).

Surprisingly enough, although there are many examples of birth-death processes on $\mathbb{N}_0$ where the spectral measure and the corresponding orthogonal polynomials are given (see \cite[Chapter 3]{MDIB} for a wide collection of examples), there is \emph{only one} explicit example of bilateral birth-death process, as far as the author knows, where the spectral matrix and the corresponding orthogonal polynomials have been explicitly computed (see \cite{ILMV}). The purpose of this paper is to compute the spectral matrix and the corresponding orthogonal polynomials of several new examples of bilateral birth-death processes and use the Karlin-McGregor representation formula to study some probabilistic properties like recurrence, the invariant distribution (if it exists) and the probability current associated with the processes.

First, in Section \ref{sec2}, we recall some results the spectral analysis of birth-death processes and bilateral birth-death processes. For birth-death processes on $\mathbb{N}_0$ we also study the absorbing $M/M/1$ queue (with constant transition rates), which will play an important role in our examples. For bilateral birth-death processes we will derive Stieltjes transform relations between the spectral matrix and the spectral measures associated with the two birth-death processes on $\mathbb{N}_0$ corresponding to the two directions to infinity. These relations will be the main tool to compute the spectral matrix in our examples. We also recall the example studied in \cite{ILMV} (also with constant transition rates) and study the \emph{symmetric} bilateral birth-death process with constant transition rates motivated by the discrete-time random walk on $\mathbb{Z}$ with an attractive or repulsive force studied in \cite[Section 6]{G1} (see also \cite{dIJ1}).

In Section \ref{sec3} we consider two cases of bilateral birth-death processes with alternating constant rates. In the first case the process will be characterized by a constant transition rate from even states and another transition rate (usually different) from odd states (see \cite{ dCIM}). The second case is similar but now the parity behavior of the birth rates will be different from the parity of the death rates. In Section \ref{sec4} we study a couple of variants of the bilateral birth-death process with constant rates studied in Section \ref{sec2}, allowing one defect at the state 0. This small variation will change the spectral analysis considerably as we will see. Finally, in Section \ref{sec5}, we will study the case where the bilateral birth-death process splits into two different absorbing $M/M/1$ queues, one in the direction to $+\infty$ and the other (with different rates) in the direction to $-\infty$. This is the most elaborated example since the spectral matrix will depend on the location of the spectrum of these independent $M/M/1$ queues.

\section{Spectral analysis of birth-death processes}\label{sec2}

In this section we recall some results concerning the spectral analysis of birth-death processes, either on $\mathbb{N}_0$ or $\mathbb{Z}$. We also recall some examples from the literature (with constant transition rates) that will be relevant in the subsequent sections.

\subsection{State space $\mathbb{N}_0$}

Let $\{X_t : t\geq0\}$ be a birth-death process on $\mathbb{N}_0$ with infinitesimal operator $\mathcal{A}$ given by the semi-infinite tridiagonal matrix
\begin{equation}\label{3mtp}
\mathcal{A}=\begin{pmatrix}
-(\lambda_0+\mu_0)&\lambda_0&&&&\\
\mu_1&-(\lambda_1+\mu_1)&\lambda_1&&&\\
&\mu_2&-(\lambda_2+\mu_2)&\lambda_2&&\\
&&\ddots&\ddots&\ddots
\end{pmatrix}.
\end{equation}
A diagram of the transitions between states is
\vspace{0.2cm}
\begin{center}
$$\begin{psmatrix}[colsep=1.9cm]
  \cnode{.4}{0}& \cnode{.4}{1} & \cnode{.4}{2}& \cnode{.4}{3}& \cnode{.4}{4}& \cnode{.4}{5} &  \rnode{6}{\Huge{\cdots}} \\
\psset{nodesep=3pt,arcangle=15,labelsep=2ex,linewidth=0.3mm,arrows=->,arrowsize=1mm
3}
\uput[u](0.95,1.95){\lambda_0}\uput[u](2.85,1.95){\lambda_1}\uput[u](4.75,1.95){\lambda_2}\uput[u](6.65,1.95){\lambda_3}\uput[u](8.55,1.95){\lambda_4}\uput[u](10.65,1.95){\lambda_5}
\uput[u](0.95,1.1){\mu_1}\uput[u](2.85,1.1){\mu_2}\uput[u](4.75,1.1){\mu_3}\uput[u](6.65,1.1){\mu_4}\uput[u](8.55,1.1){\mu_5}\uput[u](10.65,1.1){\mu_6}
 \ncarc{0}{1}\ncarc{1}{0} \ncarc{1}{2} \ncarc{2}{1}\ncarc{2}{3} \ncarc{3}{2}
\ncarc{3}{4} \ncarc{4}{3} \ncarc{4}{5}\ncarc{5}{4} \ncarc{5}{6} \ncarc{6}{5}
\psset{labelsep=-3.40ex}\nput{90}{0}{0}
\psset{labelsep=-3.40ex}\nput{90}{1}{1}
\psset{labelsep=-3.40ex}\nput{90}{2}{2}
\psset{labelsep=-3.40ex}\nput{90}{3}{3}
\psset{labelsep=-3.40ex}\nput{90}{4}{4}
\psset{labelsep=-3.40ex}\nput{90}{5}{5}
\end{psmatrix}
$$
\end{center}
\vspace{-1.2cm}
We will assume that the set of rates $\{\lambda_n, \mu_n\}$ uniquely determines the birth-death process. Define the \emph{potential coefficients} $(\pi_n)_{n\in\mathbb{N}_0}$ as usual as
\begin{equation*}\label{3potcoefN}
\pi_0=1,\quad\pi_n=\frac{\lambda_0\lambda_1\cdots \lambda_{n-1}}{\mu_1\mu_2\cdots \mu_n},\quad n\geq1.
\end{equation*}
If we assume that $\mathcal{A}$ is a closed, symmetric, self-adjoint and negative operator in the Hilbert space $\ell^2_{\pi}(\mathbb{N}_0)$ then, applying the spectral theorem (see \cite{KMc2}), we can obtain an integral representation of the  transition probability functions $P_{ij}(t)=\mathbb{P}\left(X_t=j\,|\, X_0=i\right)$ in terms of a nonnegative measure $\psi(x)$ supported on $[0,\infty)$ and a family of polynomials $(Q_n)_{n\in\mathbb{N}_0}$, generated by the three-term recurrence relation with initial conditions
\begin{equation}\label{pols1}
\begin{split}
Q_0(x)&=1,\quad Q_{-1}(x)=0,\\
-xQ_n(x)&=\lambda_nQ_{n+1}(x)-(\lambda_n+\mu_n)Q_n(x)+\mu_nQ_{n-1}(x),\quad n\in\mathbb{N}_0.
\end{split}
\end{equation}
Observe that if we define $Q(x)=(Q_0(x),Q_1(x),\ldots)^T$, then $Q(x)$ is just the eigenvector in the eigenvalue equation $-xQ(x)=\mathcal{A}Q(x)$. This integral representation, called the \emph{Karlin-McGregor formula} (see \cite{KMc2}), is given by
\begin{equation}\label{3FKMcGN}
P_{ij}(t)=\pi_j\int_{0}^{\infty}e^{-xt}Q_i(x)Q_j(x)d\psi(x),\quad i,j\in\mathbb{N}_0.
\end{equation}
In \cite{KMc3} several probabilistic properties of the birth-death processes were studied in terms of the spectral measure. For instance, the birth-death process is recurrent if and only if $\int_0^\infty x^{-1}d\psi(x)=\infty$ (also equivalent to $\sum_{n=0}^\infty (\lambda_n\pi_n)^{-1}=\infty$). Otherwise it is transient. If the process is recurrent, then it is positive recurrent if and only if the spectral measure $\psi$ has a finite jump at $x = 0$ of size $\psi(\{0\})=\left(\sum_{n=0}^\infty\pi_n\right)^{-1}$. Otherwise it is null recurrent. Other quantities like the moments of the first-passage time distributions or limit theorems can also be studied using the Karlin-McGregor representation. The reader is invited to consult \cite{MDIB} for a collection of some of these results.

A fundamental role for the birth-death process is played by the function
\begin{equation}\label{probcur0}
\Omega_{j,n}(t)=\lambda_{n-1}P_{j,n-1}(t)-\mu_nP_{j,n}(t),\quad j,n\in\mathbb{N}_0,\quad\lambda_{-1}=0,
\end{equation}
describing the \emph{probability current} in the state $n$ at time $t$, given that we start at state $j$ (see \cite[p. 383]{Gil}). $\Omega_{j,n}(t)$ represents a \emph{net probability flux} from state $n-1$ to state $n$ at time $t$. Using the Karlin-McGregor formula \eqref{3FKMcGN} and the symmetry property $\pi_n\mu_n=\lambda_{n-1}\pi_{n-1}$ we can write $\Omega_{j,n}(t)$ as
\begin{equation*}\label{probcur}
\Omega_{j,n}(t)=\mu_n\pi_n\int_0^{\infty}e^{-xt}Q_j(x)\left[Q_{n-1}(x)-Q_n(x)\right]d\psi(x),\quad j,n\in\mathbb{N}_0.
\end{equation*}
If we define the \emph{dual polynomials} $(H_n)_{n\in\mathbb{N}_0}$ (see \cite{KMc2, KMc3}) by
$$
H_0(x)=\mu_0,\quad H_{n+1}(x)=\lambda_n\pi_n\left[Q_{n+1}(x)-Q_n(x)\right],\quad n\in\mathbb{N}_0,
$$ 
then we can write $\Omega_{j,n}(t)$ as
\begin{equation*}\label{probcur2}
\Omega_{j,n}(t)=-\int_0^{\infty}e^{-xt}Q_j(x)H_n(x)d\psi(x),\quad j,n\in\mathbb{N}_0.
\end{equation*}

\begin{example}\label{3Ej7}
\emph{The absorbing $M/M/1$ queue} (\cite{KMc4}). Consider the birth-death process with constant birth-death rates given by
$$
\lambda_n=\lambda,\quad \mu_n=\mu,\quad n\in\mathbb{N}_0,\quad\lambda,\mu>0.
$$ 
We allow the state $0$ to be an absorbing state since $\mu_0=\mu$. Since the infinitesimal operator $\mathcal{A}$ in \eqref{3mtp} is the same as the infinitesimal operator of the \emph{$0$-th birth-death process} (i.e. the infinitesimal operator defined from $\mathcal{A}$ eliminating the first row
and column of $\mathcal{A}$), we can apply the identity (2.5) of \cite{KMc4} to get an explicit expression of the \emph{Stieltjes transform} of the spectral measure $\psi$, given by 
\begin{equation}\label{31server1}
B(z;\psi)\coloneqq\int_0^\infty\frac{d\psi(x)}{x-z}=\frac{\lambda+\mu-z-\sqrt{(\lambda+\mu-z)^2-4\lambda\mu}}{2\lambda\mu},
\end{equation}
where the square root is taking positive for $z<0$. Using the Perron-Stieltjes inversion formula we have that the spectral measure has only an absolutely continuous part, given by
\begin{equation}\label{mm1sp}
\psi(x)=\frac{\sqrt{(x-\sigma_-)(\sigma_+-x)}}{2\pi \lambda\mu},\quad x\in[\sigma_-,\sigma_+],
\end{equation}
where $\sigma_{\pm}=(\sqrt{\lambda}\pm\sqrt{\mu})^2$. It is easy to see that the polynomials generated by the three-term recurrence relation \eqref{pols1} are given by
$$
Q_n(x)=\left(\frac{\mu}{\lambda}\right)^{n/2}U_n\left(\frac{\lambda+\mu-x}{2\sqrt{\lambda\mu}}\right),\quad n\in\mathbb{N}_0,
$$
where $(U_n)_{n\in\mathbb{N}_0}$ are the \index{Chebychev polynomials of the second kind}Chebychev polynomials of the second kind. After making the change of variables $x=\lambda+\mu-2\sqrt{\lambda\mu}\cos\theta$ and using some basic properties of Chebychev polynomials we have that the Karlin-McGregor formula \eqref{3FKMcGN} can be written as
\begin{align*}
P_{ij}(t)&=\left(\sqrt{\frac{\lambda}{\mu}}\right)^{j-i}\int_{(\sqrt{\lambda}-\sqrt{\mu})^2}^{(\sqrt{\lambda}+\sqrt{\mu})^2}e^{-xt}U_i\left(\frac{\lambda+\mu-x}{2\sqrt{\lambda\mu}}\right)U_j\left(\frac{\lambda+\mu-x}{2\sqrt{\lambda\mu}}\right)\frac{\sqrt{4\lambda\mu-(\lambda+\mu-x)^2}}{2\pi \lambda\mu}dx\\
&=\frac{2}{\pi}\left(\sqrt{\frac{\lambda}{\mu}}\right)^{j-i}\int_0^{\pi}e^{-(\lambda+\mu-2\sqrt{\lambda\mu}\cos\theta)t}U_i(\cos\theta)U_j(\cos\theta)\sqrt{1-\cos^2\theta}\sin\theta d\theta\\
&=\frac{2}{\pi}e^{-(\lambda+\mu)t}\left(\sqrt{\frac{\lambda}{\mu}}\right)^{j-i}\int_0^\pi e^{2t\sqrt{\lambda\mu}\cos\theta}\sin((i+1)\theta)\sin((j+1)\theta)d\theta\\
&=e^{-(\lambda+\mu)t}\left(\sqrt{\frac{\lambda}{\mu}}\right)^{j-i}\left[I_{i-j}(2\sqrt{\lambda\mu}t)-I_{i+j+2}(2\sqrt{\lambda\mu}t)\right],
\end{align*}
where $I_\nu(z)$ denotes the modified Bessel function of the first kind. In the last step we have used formula (2) of \cite[pp.81]{Bate2}. This last expression seems to be new as far as the author knows. 

From the spectral measure it is possible to see that $\int_0^\infty x^{-1}\psi(x)dx<\infty$ unless $\lambda=\mu$, where it diverges. Therefore, if $\lambda\neq\mu$ the process is transient. If $\lambda=\mu$ the process is null recurrent, since the measure does not have a finite jump at $x=0$. Finally, since we have an explicit expression of $P_{ij}(t)$, we have that the probability current \eqref{probcur0} is given by
\begin{align}
\label{pcmm1}\Omega_{j,n}(t)&=e^{-(\lambda+\mu)t}\left(\sqrt{\frac{\lambda}{\mu}}\right)^{n-j}\\
\nonumber&\quad\times\left(\sqrt{\lambda\mu}\left[I_{j-n+1}(2\sqrt{\lambda\mu}t)-I_{n+j+1}(2\sqrt{\lambda\mu}t)\right]-\mu\left[I_{j-n}(2\sqrt{\lambda\mu}t)-I_{j-n}(2\sqrt{\lambda\mu}t)\right]\right).
\end{align}
In Figure \ref{fig1} this probability current is plotted as a function of $n$ starting at $j=0$ for $t=3,6,9$ and for different values of the birth-death rates $\lambda,\mu$. Note that the probability current is negative in the first states since the state 0 is absorbing (especially if $\lambda<\mu$). On the other hand, if $\lambda>\mu$ the probability current is positive for large values of the states, since the boundary $+\infty$ is attracting.
\begin{figure*}[t!]
    \centering
    \begin{subfigure}[b]{0.5\textwidth}
        \centering
        \includegraphics[height=2.5in]{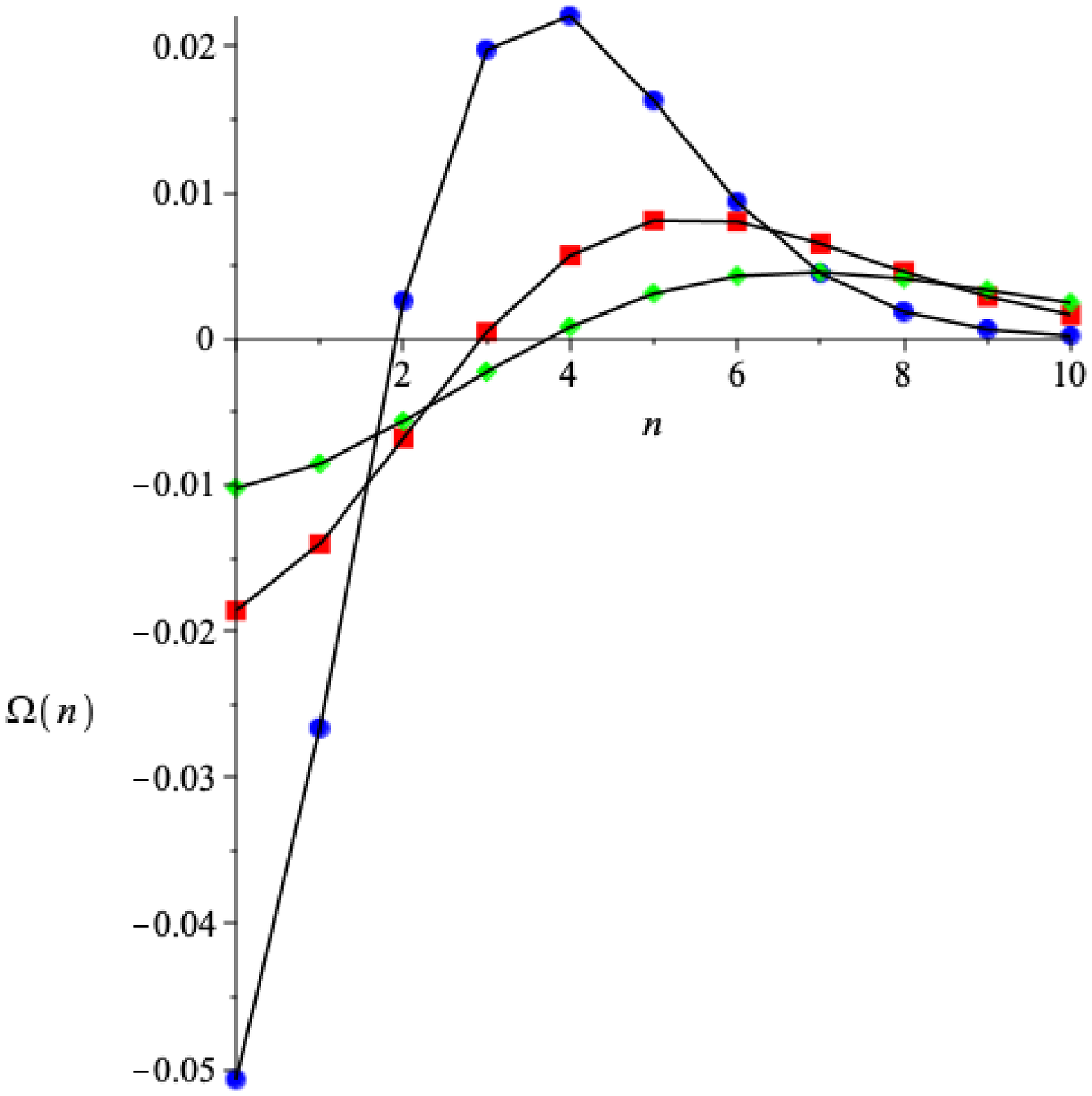}
        \caption{$\lambda=1,\mu=1$}
    \end{subfigure}%
    ~ 
    \begin{subfigure}[b]{0.5\textwidth}
        \centering
        \includegraphics[height=2.5in]{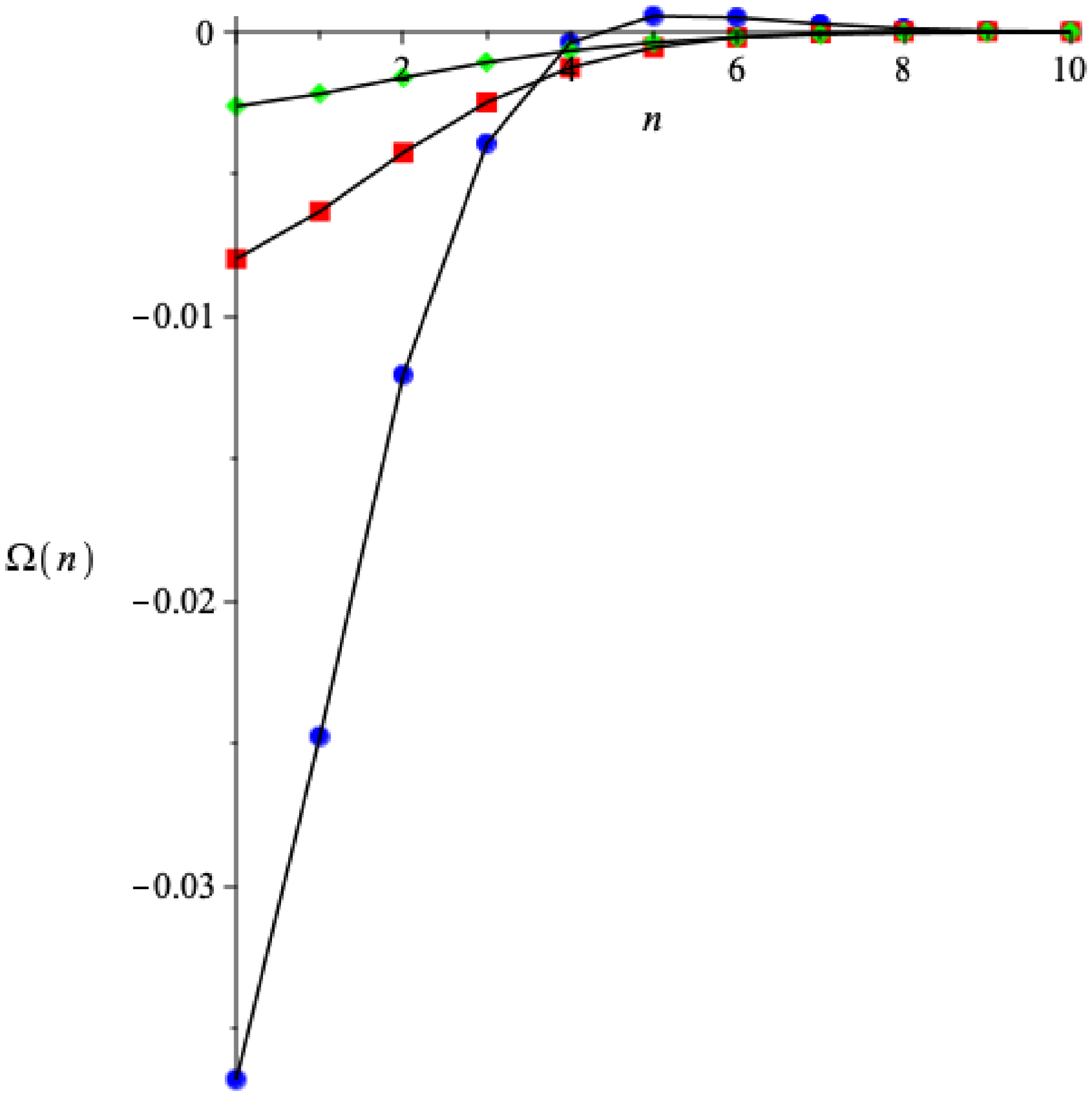}
        \caption{$\lambda=1,\mu=2$}
    \end{subfigure}\\
    \centering
    \begin{subfigure}[b]{0.5\textwidth}
        \centering
        \includegraphics[height=2.5in]{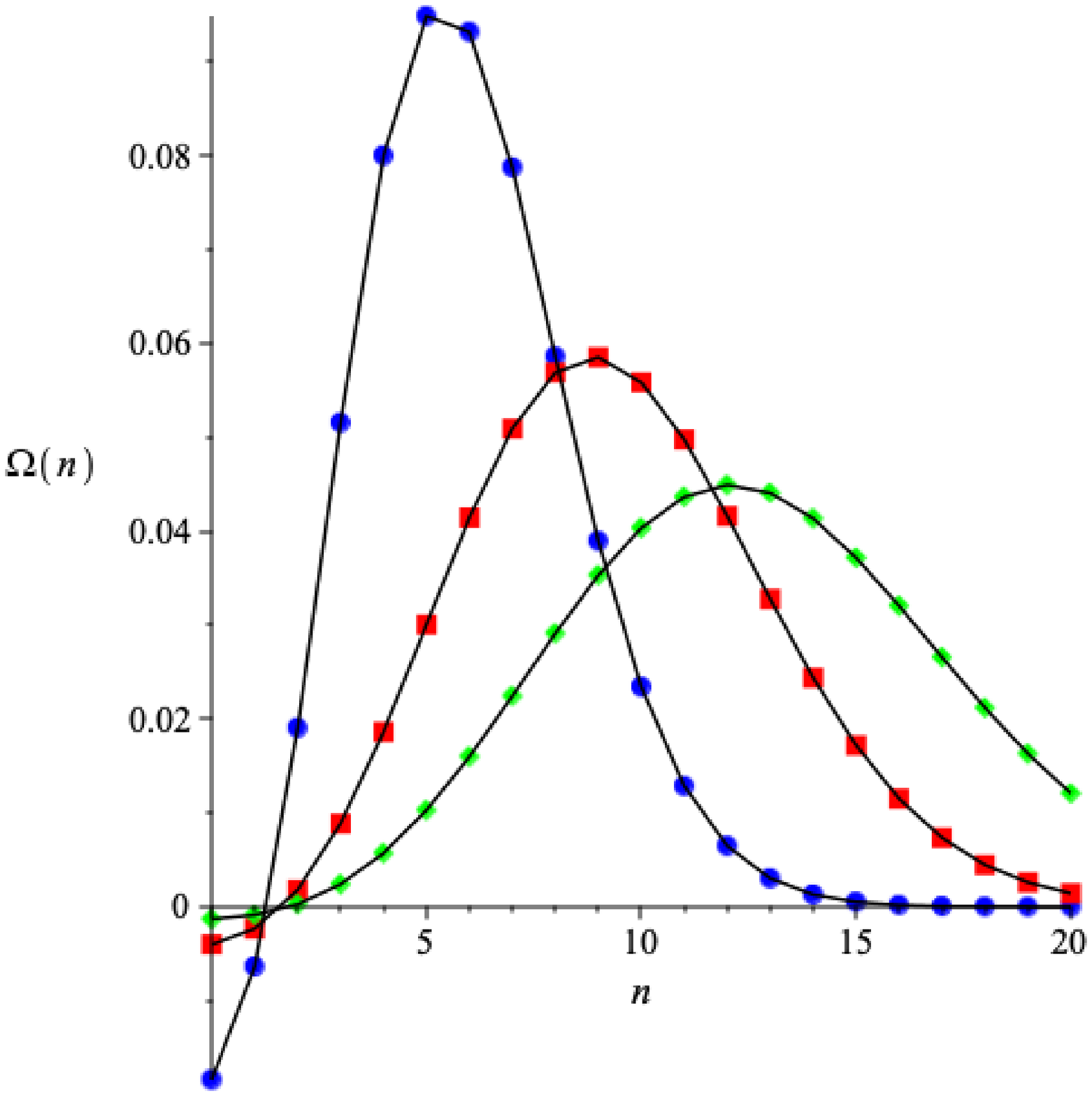}
        \caption{$\lambda=2,\mu=1$}
    \end{subfigure}%
    ~ 
    \begin{subfigure}[b]{0.5\textwidth}
        \centering
        \includegraphics[height=2.5in]{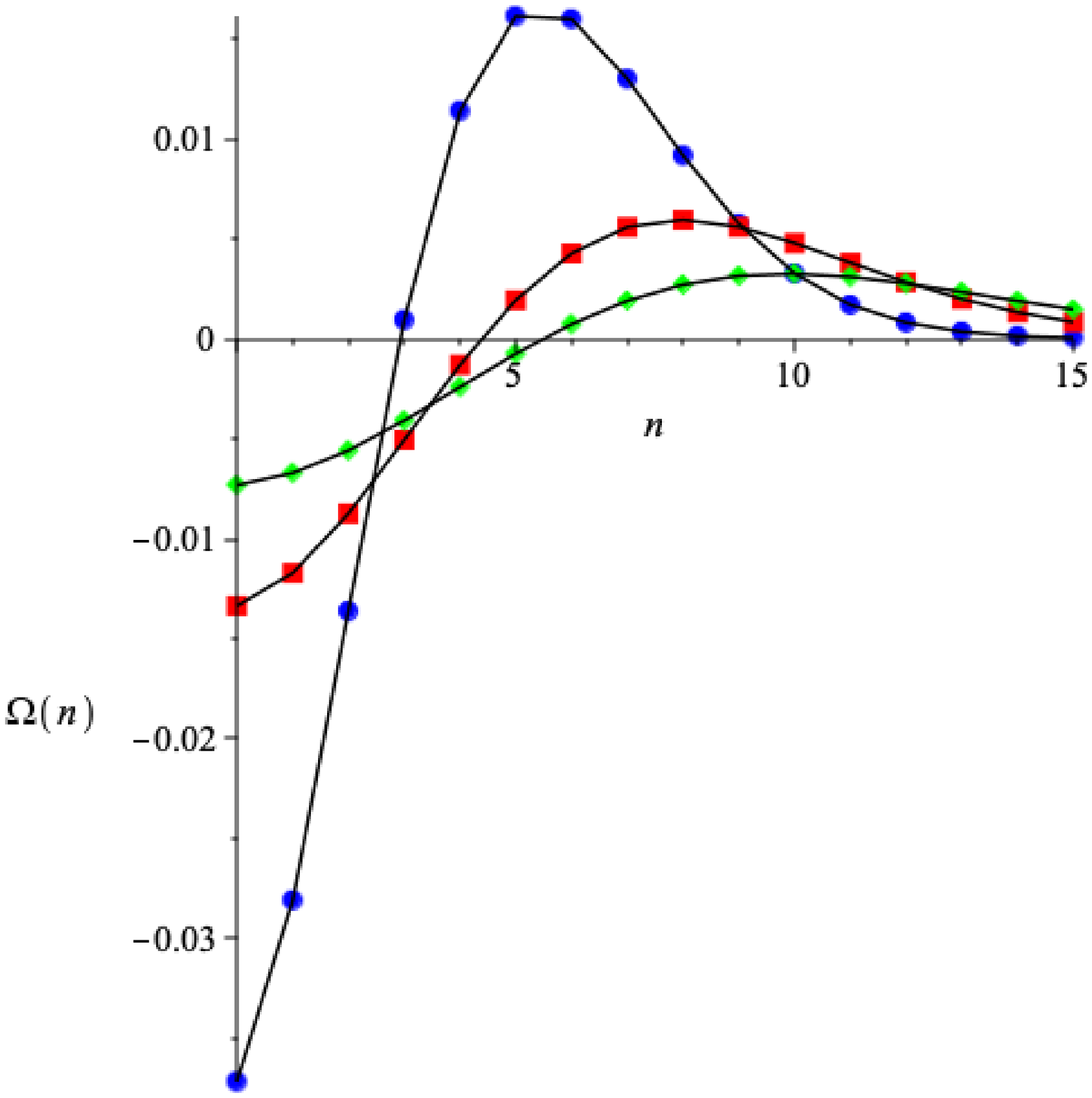}
        \caption{$\lambda=2,\mu=2$}
    \end{subfigure}
  \caption{The probability current \eqref{pcmm1} (here denoted by $\Omega(n)$) is plotted as function of $n$ for $j = 0$, $t = 3$ (blue circles), $t = 6$ (red squares) and $t =9$ (green diamonds).}  \label{fig1}
\end{figure*}

\end{example}

We invite the reader to consult \cite[Chapter 3]{MDIB} for a summary of other examples.

\subsection{State space $\mathbb{Z}$}

Let $\{X_t : t\geq0\}$ be a birth-death process on $\mathbb{Z}$ with infinitesimal operator $\mathcal{A}$ given by the doubly infinite tridiagonal matrix
\begin{equation*}\label{P1}
\mathcal{A}=\left(
\begin{array}{ccc|cccc}
\ddots&\ddots&\ddots&&&\\
&\mu_{-1}&-(\mu_{-1}+\lambda_{-1})&\lambda_{-1}&&&\\
\hline
&&\mu_{0}&-(\mu_{0}+\lambda_{0})&\lambda_{0}&&\\
&&&\mu_{1}&-(\mu_{1}+\lambda_{1})&\lambda_{1}&\\
&&&&\ddots&\ddots&\ddots
\end{array}
\right).
\end{equation*}
Now a diagram of the transitions between states is
\vspace{0.2cm}

\begin{center}
$$\begin{psmatrix}[colsep=1.9cm]
  \rnode{0}{\Huge{\cdots}}& \cnode{.4}{1} & \cnode{.4}{2}& \cnode{.4}{3}& \cnode{.4}{4}& \cnode{.4}{5} &  \rnode{6}{\Huge{\cdots}} \\
\psset{nodesep=3pt,arcangle=15,labelsep=2ex,linewidth=0.3mm,arrows=->,arrowsize=1mm
3}
\uput[u](0.95,1.95){\lambda_{-3}}\uput[u](3,1.95){\lambda_{-2}}\uput[u](4.9,1.95){\lambda_{-1}}\uput[u](6.8,1.95){\lambda_0}\uput[u](8.7,1.95){\lambda_1}\uput[u](10.8,1.95){\lambda_2}
\uput[u](0.95,1.1){\mu_{-2}}\uput[u](3,1.1){\mu_{-1}}\uput[u](4.9,1.1){\mu_0}\uput[u](6.8,1.1){\mu_1}\uput[u](8.7,1.1){\mu_2}\uput[u](10.8,1.1){\mu_3}
 \ncarc{0}{1}\ncarc{1}{0} \ncarc{1}{2} \ncarc{2}{1}\ncarc{2}{3} \ncarc{3}{2}
\ncarc{3}{4} \ncarc{4}{3} \ncarc{4}{5}\ncarc{5}{4} \ncarc{5}{6} \ncarc{6}{5}
\psset{labelsep=-3.40ex}\nput{90}{1}{-2}
\psset{labelsep=-3.40ex}\nput{90}{2}{-1}
\psset{labelsep=-3.40ex}\nput{90}{3}{0}
\psset{labelsep=-3.40ex}\nput{90}{4}{1}
\psset{labelsep=-3.40ex}\nput{90}{5}{2}
\end{psmatrix}
$$
\end{center}
\vspace{-1.3cm}
These processes are also known as \emph{bilateral birth-death processes} (following \cite{Pru}), double-ended systems or unrestricted birth-death processes (see \cite{dCIM,dCMM,dCM,GN}). Again, we will assume that the set of rates $\{\lambda_n, \mu_n\}$ uniquely determines the bilateral birth-death process. In a similar way we can define the \emph{potential coefficients} by
\begin{equation}\label{3potcoefZ}
\pi_0=1,\quad\pi_n=\frac{\lambda_0\lambda_1\cdots \lambda_{n-1}}{\mu_1\mu_2\cdots \mu_n},\quad
\pi_{-n}=\frac{\mu_0\mu_{-1}\cdots \mu_{-n+1}}{\lambda_{-1}\lambda_{-2}\cdots \lambda_{-n}},\quad n\in\mathbb{N}.
\end{equation}
If we assume that $\mathcal{A}$ is a closed, symmetric, self-adjoint and negative operator in the Hilbert space $\ell^2_{\pi}(\mathbb{Z})$ then, applying the spectral theorem (see \cite{Pru}), we can obtain an integral representation of the  transition probability functions $P_{ij}(t)=\mathbb{P}\left(X_t=j\,|\, X_0=i\right)$ in terms of \emph{three measures} $\psi_{11}(x), \psi_{22}(x)$ and $\psi_{12}(x)$ supported on $[0,\infty)$ and two linearly independent families of polynomials $(Q_n^{\alpha})_{n\in\mathbb{Z}},\alpha=1,2,$ generated by the three-term recurrence relations with initial conditions
\begin{align}
\nonumber Q_0^1(x)&=1,\quad Q_{0}^2(x)=0,\\
\label{3TTRRZ}Q_{-1}^1(x)&=0,\quad Q_{-1}^2(x)=1,\\
\nonumber -xQ_n^{\alpha}(x)&=\lambda_nQ_{n+1}^{\alpha}(x)-(\lambda_n+\mu_n)Q_n^{\alpha}(x)+\mu_nQ_{n-1}^{\alpha}(x),\quad n\in\mathbb{Z}.
\end{align}
This integral representation, called the \emph{Karlin-McGregor formula}, is given by
\begin{equation}\label{3FKMcGZ}
P_{ij}(t)=\pi_j\int_{0}^{\infty}e^{-xt}\sum_{\alpha,\beta=1}^{2}Q_i^{\alpha}(x)Q_j^{\beta}(x)d\psi_{\alpha\beta}(x),\quad i,j\in\mathbb{Z}.
\end{equation}
These three measures can be grouped in a positive definite $2\times2$ matrix, called the \emph{spectral matrix} of the bilateral birth-death process:
$$
\Psi(x)=\begin{pmatrix} \psi_{11}(x) & \psi_{12}(x)\\ \psi_{12}(x) & \psi_{22}(x) \end{pmatrix}.
$$
Therefore, the Karlin-McGregor formula \eqref{3FKMcGZ} can be written in matrix form as
\begin{equation}\label{3FKMcGZ2}
P_{ij}(t)=\pi_j\int_{0}^{\infty}e^{-xt}\left(Q_i^1(x),Q_i^2(x)\right)d\Psi(x)\begin{pmatrix}Q_j^1(x)\\Q_j^2(x)\end{pmatrix},\quad i,j\in\mathbb{Z}.
\end{equation}
The computation of the measures $\psi_{\alpha\beta}(x),\alpha,\beta=1,2$ can be reduced to study two birth-death processes on $\mathbb{N}_0$ corresponding to the two directions to infinity, with infinitesimal operators $\mathcal{A}_{ij}^+=\mathcal{A}_{ij}, i,j\geq0$ and $\mathcal{A}_{ij}^-=\mathcal{A}_{ij}, i,j\leq-1$, i.e.
\begin{equation}\label{AAp}
\mathcal{A}^+=\begin{pmatrix}
-(\lambda_0+\mu_0)&\lambda_0&&&\\
\mu_1&-(\lambda_1+\mu_1)&\lambda_1&&\\
&\mu_1&-(\lambda_2+\mu_2)&\lambda_2&\\
&&\ddots&\ddots&\ddots
\end{pmatrix},
\end{equation}
and
\begin{equation}\label{AAm}
\mathcal{A}^-=\begin{pmatrix}
-(\lambda_{-1}+\mu_{-1})&\mu_{-1}&&&\\
\lambda_{-2}&-(\lambda_{-2}+\mu_{-2})&\mu_{-2}&&\\
&\lambda_{-3}&-(\lambda_{-3}+\mu_{-3})&\mu_{-3}&\\
&&\ddots&\ddots&\ddots
\end{pmatrix}.
\end{equation}
Let us denote by  $\psi^{\pm}$ the spectral measures associated with $\mathcal{A}^{\pm}$. In the last section of \cite{KMc6} a method to relate the Stieltjes transforms of $\psi_{\alpha\beta},\alpha,\beta=1,2$ in terms of the Stieltjes transforms of  $\psi^{\pm}$ was given for discrete-time birth-death chains using probabilistic arguments. We will use here the same arguments to derive similar relations. For that, let us call 
$$
F_{ij}(t)=\mathbb{P}\left(X_{\tau}=j\,\mbox{for some $\tau$}, 0<\tau\leq t\,|\,X_0=i\right), \quad i\neq j,
$$
the first-passage time distributions and
$$
F_{ii}(t)=\mathbb{P}\left(X_{\tau_1}\neq i, X_{\tau_2}=i\,\mbox{for some $\tau_1,\tau_2$}, 0<\tau_1<\tau_2\leq t\,|\,X_0=i\right),
$$
the recurrence time distributions. Let us denote $\widehat{P}_{ij}(s)$ and $\widehat{F}_{ij}(s)$ the Laplace transforms of $P_{ij}(t)$ and $F_{ij}(t)$, respectively, i.e.
\begin{equation}\label{lapPF}
\widehat{P}_{ij}(s)=\int_0^\infty e^{-st}P_{ij}(t)dt,\quad\widehat{F}_{ij}(s)=\int_0^\infty e^{-st}dF_{ij}(t),\quad i,j\in\mathbb{Z}.
\end{equation}
We will use the same notation for the birth-death processes on $\mathbb{N}_0$ generated by $\mathcal{A}^{\pm}$ (i.e. $P_{ij}^{\pm}(t),F_{ij}^{\pm}(t)$ and $\widehat{P}_{ij}^{\pm}(s),\widehat{F}_{ij}^{\pm}(s)$). The Laplace transforms $\widehat{P}_{ij}(s)$ and $\widehat{F}_{ij}(s)$ are related by the following formulas (see \cite{KMc2,KMc3})
\begin{equation}\label{reneweq}
\begin{split}
\widehat{P}_{ii}(s)&=\frac{1}{s+\lambda_i+\mu_i}+\widehat{P}_{ii}(s)\widehat{F}_{ij}(s),\\
\widehat{P}_{ij}(s)&=\widehat{P}_{jj}(s)\widehat{F}_{ij}(s),\quad i\neq j.
\end{split}
\end{equation}
From the identities
\begin{align*}
\widehat{F}_{00}(s)&=\widehat{F}_{00}^+(s)+\frac{\mu_0}{s+\lambda_0+\mu_0}\widehat{F}_{-1,0}(s),\\
\widehat{F}_{-1,0}(s)&=\lambda_{-1}\widehat{P}_{-1,-1}^-(s),\\
\widehat{F}_{00}^+(s)&=1-\frac{1}{(s+\lambda_0+\mu_0)\widehat{P}_{00}^+(s)},
\end{align*}
it is found that
\begin{equation}\label{form11}
\widehat{P}_{00}(s)=\frac{\widehat{P}_{00}^+(s)}{1-\lambda_{-1}\mu_0\widehat{P}_{00}^+(s)\widehat{P}_{-1,-1}^-(s)}.
\end{equation}
Similarly, from the identities
\begin{align*}
\widehat{F}_{-1,-1}(s)&=\widehat{F}_{-1,-1}^-(s)+\frac{\lambda_{-1}}{s+\lambda_{-1}+\mu_{-1}}\widehat{F}_{0,-1}(s),\\
\widehat{F}_{0,-1}(s)&=\mu_{0}\widehat{P}_{00}^+(s),\\
\widehat{F}_{-1,-1}^-(s)&=1-\frac{1}{(s+\lambda_{-1}+\mu_{-1})\widehat{P}_{-1,-1}^-(s)},
\end{align*}
we obtain
\begin{equation}\label{form22}
\widehat{P}_{-1,-1}(s)=\frac{\widehat{P}_{-1,-1}^-(s)}{1-\lambda_{-1}\mu_0\widehat{P}_{00}^+(s)\widehat{P}_{-1,-1}^-(s)}.
\end{equation}
Finally, $\widehat{P}_{-1,0}(s)=\widehat{P}_{00}(s)\widehat{F}_{-1,0}(s)$ gives
\begin{equation}\label{form12}
\widehat{P}_{-1,0}(s)=\frac{\lambda_{-1}\widehat{P}_{00}^+(s)\widehat{P}_{-1,-1}^-(s)}{1-\lambda_{-1}\mu_0\widehat{P}_{00}^+(s)\widehat{P}_{-1,-1}^-(s)}.
\end{equation}
Now, using the Karlin-McGregor formula \eqref{3FKMcGZ} in \eqref{lapPF} we obtain
\begin{equation}\label{pijsb}
\widehat{P}_{ij}(s)=\pi_j\int_{0}^{\infty}\frac{1}{x+s}\sum_{\alpha,\beta=1}^{2}Q_i^{\alpha}(x)Q_j^{\beta}(x)d\psi_{\alpha\beta}(x)=\pi_jB\left(-s;\sum_{\alpha,\beta=1}^{2}Q_i^{\alpha}Q_j^{\beta}\psi_{\alpha\beta}\right),
\end{equation}
where $B(z;\psi)$ is defined by $B(z;\psi)=\int_0^\infty(x-z)^{-1}d\psi(x)$. Similar formulas hold for $\widehat{P}_{ij}^{\pm}(s)$ with the measures $\psi^{\pm}$. Therefore formulas \eqref{form11}, \eqref{form22} and \eqref{form12} are equivalent to the Stieltjes transform relations
\begin{align}
\nonumber B(z;\psi_{11})&=\frac{B(z;\psi^+)}{1-\lambda_{-1}\mu_0B(z;\psi^+)B(z;\psi^-)},\\
\label{3BzZ}\frac{\mu_0}{\lambda_{-1}}B(z;\psi_{22})&=\frac{B(z;\psi^-)}{1-\lambda_{-1}\mu_0B(z;\psi^+)B(z;\psi^-)},\\
\nonumber B(z;\psi_{12})&=\frac{\lambda_{-1}B(z;\psi^+)B(z;\psi^-)}{1-\lambda_{-1}\mu_0B(z;\psi^+)B(z;\psi^-)}.
\end{align}
Relations \eqref{form11}, \eqref{form22} and \eqref{form12} can also be derived from the tools developed in \cite{GV1}. Different arguments to compute the spectral measures $\psi_{\alpha\beta}(x),\alpha,\beta=1,2$ were given in \cite{Pru}, using the asymptotic analysis of the corresponding orthogonal polynomials, or in \cite{ILMV}, using tools from the spectral theory of self-adjoint operators.

\begin{remark}
Observe that the equations \eqref{reneweq} give a direct relation between the Laplace transforms of $P_{ij}(t)$ and the first-passage time distributions $F_{ij}(t)$. In \cite[pp. 64]{PruT} (see also \cite[Theorem 3.2]{Pru}) one can find an explicit formula for $\widehat{P}_{ij}(s)$ in terms of the two families of polynomials $(Q_n^\alpha)_{n\in\mathbb{Z}},\alpha=1,2$ defined by \eqref{3TTRRZ}, the potential coefficients $(\pi_n)_{n\in\mathbb{Z}}$ defined by \eqref{3potcoefZ} and $\widehat{P}_{00}(s), \widehat{P}_{0,-1}(s),\widehat{P}_{-1,0}(s)$ and $\widehat{P}_{-1,-1}(s)$. Since these last four expressions can be written in terms of the Stieltjes transforms $B(z;\psi_{\alpha\beta}),\alpha,\beta=1,2$ (see \eqref{pijsb}), if we have explicit formulas for the Stieltjes transforms, then we will have explicit formulas for $\widehat{P}_{ij}(s)$ and consequently for $\widehat{F}_{ij}(s)$. This gives an alternative way of computing the first-passage time distributions $F_{ij}(t)$ by applying the inverse Laplace transform.
\end{remark}

As in the case of birth-death processes on $\mathbb{N}_0$ we can derive some probabilistic properties of bilateral birth-death processes in terms of the spectral matrix. This was partially done in \cite{PruT} for recurrence, limit theorems or absorption probabilities, but for some reason it did not appear in \cite{Pru}. Recently, since a bilateral birth-death processes is a special case of a \emph{quasi-birth-and-death process} with state space $\mathbb{N}_0\times\{1,2\}$ (see \cite{LaR} for more information about these processes), some other probabilistic properties were derived in \cite{DR}. In particular, using Corollary 4.7 of \cite{DR} for $\alpha=0$, we have that the process is recurrent if and only if 
\begin{equation}\label{rec}
\int_{0}^\infty\frac{\psi_{11}(x)}{x}dx=\infty,\quad\mbox{or}\quad\pi_{-1}\int_{0}^\infty\frac{\psi_{22}(x)}{x}dx=\infty, 
\end{equation}
and positive recurrent if and only if 
\begin{equation}\label{recp}
\mbox{either}\quad \psi_{11}(x)\quad\mbox{or}\quad \pi_{-1}\psi_{22}(x)\quad \mbox{has a jump at the point 0}.
\end{equation}
The size of this jump is the same for all measures and given by $\psi(\{0\})=\left(\sum_{n\in\mathbb{Z}}\pi_n\right)^{-1}$ (see \cite[p. 119]{PruT}). Other quantities like the first moment of the first-passage time distribution or limit theorems can also be studied using the spectral matrix (see \cite[Ch. 5]{PruT} for more information).

Again, a fundamental role for bilateral birth-death processes is played by the \emph{probability current}
\begin{equation}\label{probcur0Z}
\Omega_{j,n}(t)=\lambda_{n-1}P_{j,n-1}(t)-\mu_nP_{j,n}(t),\quad j,n\in\mathbb{Z}.
\end{equation}
Using the Karlin-McGregor formula \eqref{3FKMcGZ} and the symmetry property $\pi_n\mu_n=\lambda_{n-1}\pi_{n-1}$ we can write $\Omega_{j,n}(t)$ as
\begin{equation}\label{probcurZ}
\Omega_{j,n}(t)=\mu_n\pi_n\int_0^{\infty}e^{-xt}\sum_{\alpha,\beta=1}^2Q_j\alpha(x)\left[Q_{n-1}^\beta(x)-Q_n^\beta(x)\right]d\psi_{\alpha\beta}(x),\quad j,n\in\mathbb{Z}.
\end{equation}
Again, if we define the \emph{dual polynomials} $(H_n^\alpha)_{n\in\mathbb{Z}},\alpha=1,2,$ (see \cite{PruT,Pru}) by
$$
H_{n+1}^\alpha(x)=\lambda_n\pi_n\left[Q_{n+1}^\alpha(x)-Q_n^\alpha(x)\right],\quad \alpha=0,1,\quad n\in\mathbb{Z},
$$ 
then we can write $\Omega_{j,n}(t)$ as
\begin{equation*}\label{probcur2Z}
\Omega_{j,n}(t)=-\int_0^{\infty}e^{-xt}\sum_{\alpha,\beta=1}^2Q_j^\alpha(x)H_n^\beta(x)d\psi_{\alpha\beta}(x),\quad j,n\in\mathbb{Z}.
\end{equation*}

As we pointed out in the Introduction, there is \emph{only one} explicit example of bilateral birth-death process, as far as the author knows, where the spectral matrix and the corresponding orthogonal polynomials have been explicitly computed (see \cite{ILMV}). We will recall that example here and we will also give another simple example motivated by the discrete-time random walk on $\mathbb{Z}$ with an attractive or repulsive force studied in \cite[Section 6]{G1} (see also \cite{dIJ1}).

\begin{example}\label{exb1}
\emph{Bilateral birth-death process with constant rates} (\cite{ILMV}). Consider the bilateral birth-death process with constant birth-death rates  
$$
\lambda_n=\lambda,\quad\mu_n=\mu,\quad n\in\mathbb{Z},\quad\lambda,\mu>0.
$$
The matrix $\mathcal{A}^+$ in \eqref{AAp} is the same as the one of the absorbing $M/M/1$ queue in Example \ref{3Ej7}, while $\mathcal{A}^-$ in \eqref{AAm} is the symmetric matrix of $\mathcal{A}^+$. Therefore both processes generate the same Stieltjes transform given by \eqref{31server1}. Following \eqref{3BzZ} and rationalizing we obtain that
\begin{align*}
B(z;\psi_{11})&=\frac{\mu}{\lambda}B(z;\psi_{22})=\frac{1}{\sqrt{(\lambda+\mu-z)^2-4\lambda\mu}},\\
B(z;\psi_{12})&=\frac{1}{2\mu}\left(-1+\frac{\lambda+\mu-z}{\sqrt{(\lambda+\mu-z)^2-4\lambda\mu}}\right).
\end{align*}
Observe that the jumps, if any, should be located at $x=\sigma_{\pm}$ where $\sigma_{\pm}=(\sqrt{\lambda}\pm\sqrt{\mu})^2$. However it follows easily that the size of these jumps must be 0, so there are no jumps. The spectral matrix is then given by
\begin{equation*}
\Psi(x)=\frac{1}{\pi\sqrt{(x-\sigma_-)(\sigma_+-x)}}\begin{pmatrix}1&(\lambda+\mu-x)/2\mu\\(\lambda+\mu-x)/2\mu&\lambda/\mu\end{pmatrix},\quad x\in[\sigma_-,\sigma_+].
\end{equation*}
The polynomials generated by the three-term recurrence relation \eqref{3TTRRZ} (something that it was not pointed out in \cite{ILMV}) are given by
\begin{align*}
Q_n^1(x)&=\left(\frac{\mu}{\lambda}\right)^{n/2}U_n\left(\frac{\lambda+\mu-x}{2\sqrt{\lambda\mu}}\right),\quad Q_{-n-1}^1(x)=-\left(\frac{\lambda}{\mu}\right)^{(n+1)/2}U_{n-1}\left(\frac{\lambda+\mu-x}{2\sqrt{\lambda\mu}}\right),\quad n\in\mathbb{N}_0,\\
Q_n^2(x)&=-\left(\frac{\mu}{\lambda}\right)^{(n+1)/2}U_{n-1}\left(\frac{\lambda+\mu-x}{2\sqrt{\lambda\mu}}\right),\quad Q_{-n-1}^2(x)=\left(\frac{\lambda}{\mu}\right)^{n/2}U_{n}\left(\frac{\lambda+\mu-x}{2\sqrt{\lambda\mu}}\right),\quad n\in\mathbb{N}_0,
\end{align*}
where again $(U_n)_{n\in\mathbb{N}_0}$ are the Chebychev polynomials of the second kind. The transition probability functions $P_{ij}(t)$ can then be approximated using the Karlin-McGregor formula \eqref{3FKMcGZ}. In this case the transition probability functions were explicitly computed in \cite{Con}, given by
$$
P_{ij}(t)=e^{-(\lambda+\mu)t}\left(\sqrt{\frac{\lambda}{\mu}}\right)^{j-i}I_{j-i}\left(2\sqrt{\lambda\mu}t\right),\quad i,j\in\mathbb{Z},
$$
where again $I_\nu(z)$ denotes the modified Bessel function of the first kind. The previous formula can also be derived using basic properties of Chebychev polynomials in the Karlin-McGregor formula \eqref{3FKMcGZ}.

From the spectral matrix it is possible to see that $\int_0^\infty x^{-1}\psi_{11}(x)dx<\infty$ and $\int_0^\infty x^{-1}\psi_{22}(x)dx<\infty$ unless $\lambda=\mu$, where both integrals diverge. Therefore, if $\lambda\neq\mu$ the process is transient. If $\lambda=\mu$ the process is null recurrent, since the measure does not have a jump at $x=0$. Finally, since we have an explicit expression of $P_{ij}(t)$, we have that the probability current \eqref{probcur0Z} is given by
$$
\Omega_{j,n}(t)=\mu e^{-(\lambda+\mu)t}\left(\sqrt{\frac{\lambda}{\mu}}\right)^{n-j}\left(\sqrt{\lambda/\mu}I_{n-j-1}(2\sqrt{\lambda\mu}t)-I_{n-j}(2\sqrt{\lambda\mu}t)\right),\quad j,n\in\mathbb{Z}.
$$
Some plots of $\Omega_{0,n}(t)$ can be found in Figure 9 of \cite{GN}.

\end{example}

\begin{example}\label{exb2}
\emph{Symmetric bilateral birth-death process with constant rates}. Consider the bilateral birth-death process with birth-death rates  
$$
\lambda_n=\lambda,\quad\mu_n=\mu, \quad n\in\mathbb{N}_0,\quad \lambda_{-n}=\mu,\quad \mu_{-n}=\lambda, \quad n\in\mathbb{N},\quad \lambda,\mu>0. 
$$
The matrices $\mathcal{A}^{\pm}$ in \eqref{AAp} and \eqref{AAm} are equal and the same as the one of the absorbing $M/M/1$ queue in Example \ref{3Ej7}. Therefore the corresponding Stieltjes transforms are given by \eqref{31server1}. Following \eqref{3BzZ} and rationalizing we obtain that
\begin{align*}
B(z;\psi_{11})&=B(z;\psi_{22})=\frac{(\lambda-\mu)(\lambda+\mu-z)-(\lambda+\mu)\sqrt{(\lambda+\mu-z)^2-4\lambda\mu}}{2\mu z(2\lambda+2\mu-z)},\\
B(z;\psi_{12})&=\frac{\mu(\mu+2z)-(\lambda-z)^2+(\lambda+\mu-z)\sqrt{(\lambda+\mu-z)^2-4\lambda\mu}}{2\mu z(2\lambda+2\mu-z)}.
\end{align*}
The spectral matrix $\Psi(x)=\Psi_c(x)+\Psi_d(x)$ has now an absolutely continuous part $\Psi_c(x)$, given by 
\begin{equation*}
\Psi_c(x)=\frac{\sqrt{(x-\sigma_-)(\sigma_+-x)}}{2\pi\mu x(2\lambda+2\mu-x)}\begin{pmatrix}\lambda+\mu&\lambda+\mu-x\\\lambda+\mu-x&\lambda+\mu\end{pmatrix},\quad x\in[\sigma_-,\sigma_+],
\end{equation*}
where $\sigma_{\pm}=(\sqrt{\lambda}\pm\sqrt{\mu})^2$ and a discrete part $\Psi_d(x)$, given by 
\begin{equation*}
\Psi_d(x)=\frac{\mu-\lambda}{2\mu}\left[\begin{pmatrix}1&1\\1&1\end{pmatrix}\delta_0(x)+\begin{pmatrix}1&-1\\-1&1\end{pmatrix}\delta_{2\lambda+2\mu}(x)\right]\mathbf{1}_{\{\mu>\lambda\}},
\end{equation*}
where $\mathbf{1}_A$ is the indicator function and $\delta_a(x)$ is the Dirac delta located at $x=a$. The polynomials generated by the three-term recurrence relation \eqref{3TTRRZ} are given by
\begin{align*}
Q_n^1(x)&=\left(\frac{\mu}{\lambda}\right)^{n/2}U_n\left(\frac{\lambda+\mu-x}{2\sqrt{\lambda\mu}}\right),\quad Q_{-n-1}^1(x)=-\left(\frac{\mu}{\lambda}\right)^{(n+1)/2}U_{n-1}\left(\frac{\lambda+\mu-x}{2\sqrt{\lambda\mu}}\right),\quad n\in\mathbb{N}_0,\\
Q_n^2(x)&=-\left(\frac{\mu}{\lambda}\right)^{(n+1)/2}U_{n-1}\left(\frac{\lambda+\mu-x}{2\sqrt{\lambda\mu}}\right),\quad Q_{-n-1}^2(x)=\left(\frac{\mu}{\lambda}\right)^{n/2}U_{n}\left(\frac{\lambda+\mu-x}{2\sqrt{\lambda\mu}}\right),\quad n\in\mathbb{N}_0,
\end{align*}
where again $(U_n)_{n\in\mathbb{N}_0}$ are the Chebychev polynomials of the second kind. The transition probability functions $P_{ij}(t)$ can then be approximated using the Karlin-McGregor formula \eqref{3FKMcGZ}. Unlike the previous example, no explicit formula for $P_{ij}(t)$ has been found in terms of Bessel functions, as far as the author knows.
\begin{remark}
Observe that this example is different from the one studied in Section 5.3 of \cite{GN}, where $\mu_0=\lambda$ (here $\mu_0=\mu$). The rest of birth-death rates are the same. This small variation modifies the spectral matrix. We will get back to this variation later at the end of Section \ref{sec4}.
\end{remark}
Following \eqref{rec}, from the spectral matrix it is possible to see that if $\lambda>\mu$ then $\int_0^\infty x^{-1}\psi_{11}(x)dx<\infty$ and $\int_0^\infty x^{-1}\psi_{22}(x)dx<\infty$. This is because the point $x=0$ never belongs to the support of the spectral matrix. Therefore, if $\lambda>\mu$ the process is transient. Otherwise, if $\lambda\leq\mu$, both integrals diverge because the point $x=0$ belongs to the support of the spectral matrix, either if we have a discrete Dirac delta at $x=0$ (for $\lambda<\mu$) or the absolutely continuous support reaches $x=0$ (for $\lambda=\mu$). Therefore, if $\lambda\leq\mu$ the process is recurrent. Following \eqref{recp}, if $\lambda<\mu$, both measures $\psi_{11},\psi_{22}$ will always have a jump at the point 0. Therefore it will be positive recurrent. If $\lambda=\mu$ then the process will be null recurrent. Since the potential coefficients are given here by (see \eqref{3potcoefZ})
$$
\pi_n=\left(\frac{\lambda}{\mu}\right)^{n},\quad \pi_{-n-1}=\left(\frac{\lambda}{\mu}\right)^{n},\quad n\in\mathbb{N}_0,
$$
we have that the invariant distribution $\bm\pi$ for this process is given by
$$
\bm\pi=\frac{\mu-\lambda}{2\mu}\left(\cdots,\frac{\lambda}{\mu},1,1,\frac{\lambda}{\mu},\frac{\lambda^2}{\mu^2},\cdots\right),\quad \mu>\lambda.
$$

Now we do not have an explicit expression of the transition probability functions, so we also do not have an explicit expression of the probability current \eqref{probcur0Z}. Nevertheless we can get an approximation using \eqref{probcurZ}, the spectral matrix and the corresponding orthogonal polynomials. In Figure \ref{fig2} this probability current is plotted as a function of $n$ starting at $j=0$ for $t=3,6,9$ and for different values of the birth-death rates $\lambda,\mu$. Note that for $\lambda<\mu$ the endpoints $\pm\infty$ are reflecting boundaries but if $\lambda>\mu$ the probability current shifts to the left if $n < 0$ and to the right for $n > 0$ as $t$ increases, since $\pm\infty$ absorbing boundaries.
\begin{figure*}[t!]
    \centering
    \begin{subfigure}[b]{0.5\textwidth}
        \centering
        \includegraphics[height=2.5in]{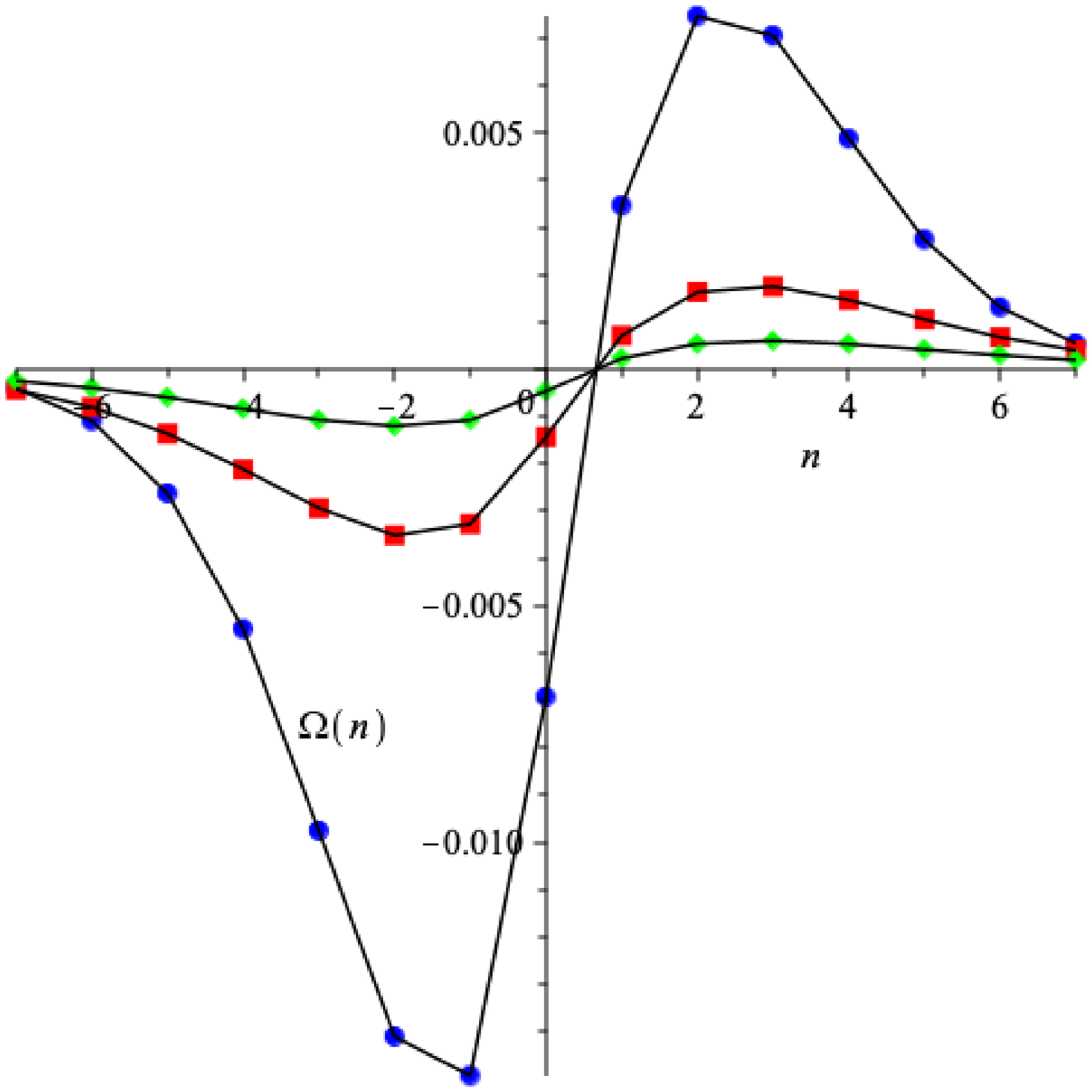}
        \caption{$\lambda=1,\mu=2$}
    \end{subfigure}%
    ~ 
    \begin{subfigure}[b]{0.5\textwidth}
        \centering
        \includegraphics[height=2.5in]{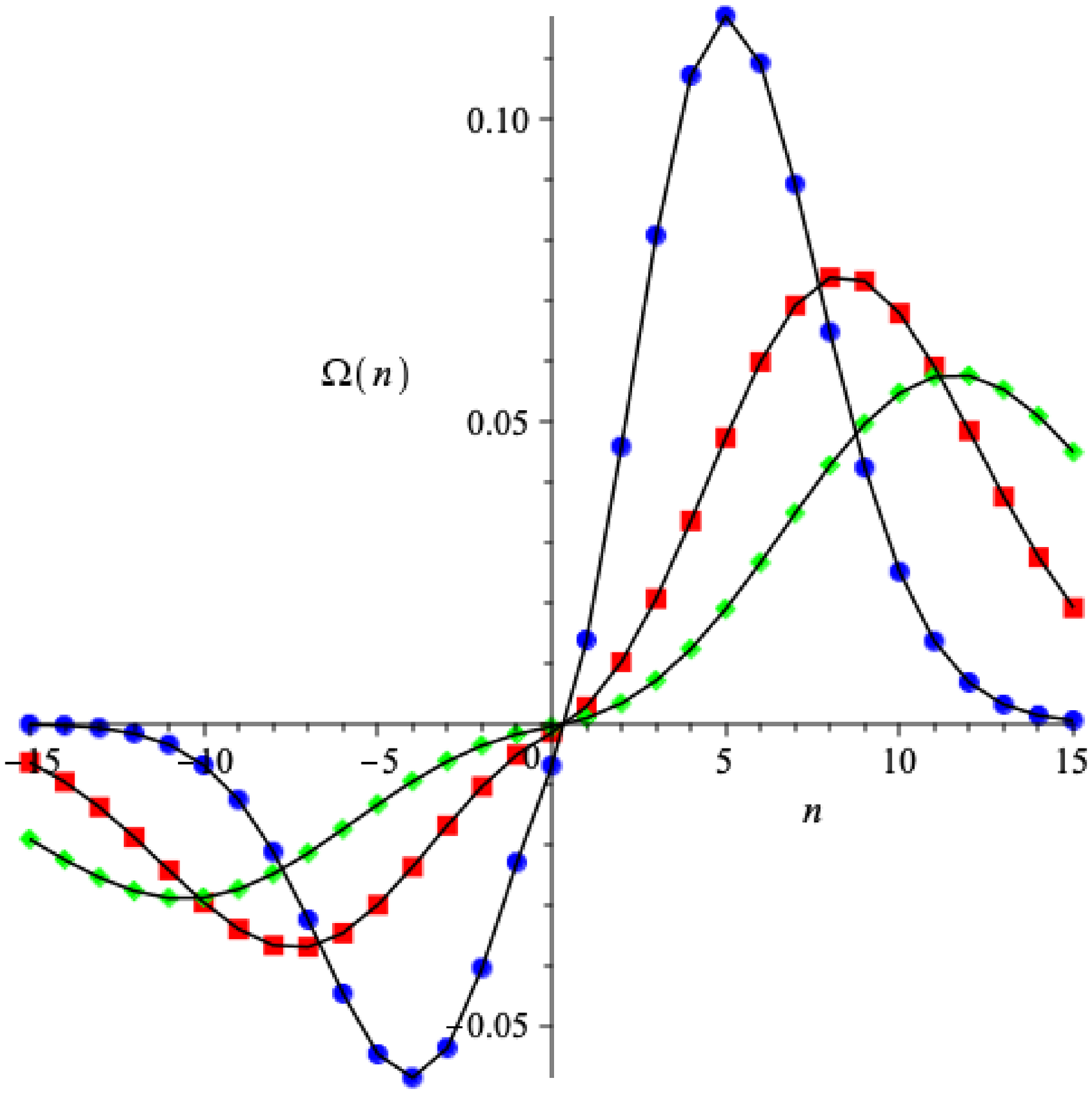}
        \caption{$\lambda=2,\mu=1$}
    \end{subfigure}
  \caption{The probability current \eqref{probcurZ} (here denoted by $\Omega(n)$) is plotted as function of $n$ for $j = 0$, $t = 3$ (blue circles), $t = 6$ (red squares) and $t =9$ (green diamonds).}  \label{fig2}
\end{figure*}

\end{example}

\section{Bilateral birth-death process with alternating constant rates}\label{sec3}

In this section we will study a couple of examples of bilateral birth-death processes with alternating constant rates. We will distinguish two cases. In the first case the process will be characterized by a constant transition rate $\lambda$ from even states and another transition rate $\mu$ from odd states (see \cite{ dCIM}). The second case is similar but now the parity behavior of the birth rates will be different from parity of the death rates. The infinitesimal operators associated with these processes are also known as Jacobi matrices with periodic recurrence coefficients (period 2 in this case) and have been extensively studied in the area of orthogonal polynomials (see for instance \cite{vAssB}).

\subsection{Case 1}
Consider the bilateral birth-death process with birth-death rates given by
$$
\lambda_{2n}=\lambda,\quad\lambda_{2n+1}=\mu, \quad \mu_{2n}=\lambda,\quad \mu_{2n+1}=\mu, \quad n\in\mathbb{Z},\quad \lambda,\mu>0. 
$$
The matrices $\mathcal{A}^{\pm}$ in \eqref{AAp} and \eqref{AAm} are now given by
\begin{equation*}
\mathcal{A}^+=\begin{pmatrix}
-2\lambda&\lambda&&&\\
\mu&-2\mu&\mu&&\\
&\lambda&-2\lambda&\lambda&\\
&&\ddots&\ddots&\ddots
\end{pmatrix},\quad
\mathcal{A}^-=\begin{pmatrix}
-2\mu&\mu&&&\\
\lambda&-2\lambda&\lambda&&\\
&\mu&-2\mu&\mu&\\
&&\ddots&\ddots&\ddots
\end{pmatrix}.
\end{equation*}
Observe that $\mathcal{A}^-$ is the infinitesimal operator of the 0-th birth-death process on $\mathbb{N}_0$ associated to $\mathcal{A}^+$ (i.e. the infinitesimal operator defined from $\mathcal{A^+}$ eliminating the first row and column of $\mathcal{A^+}$). Also $\mathcal{A}^-$ is the same matrix as $\mathcal{A}^+$ by replacing $\lambda$ by $\mu$. Applying twice the identity (2.5) of \cite{KMc4} we obtain the following algebraic relation satisfied by the Stieltjes transform of the measure $\psi^+$ associated with $\mathcal{A}^+$:
$$
\lambda\mu(z-2\lambda)B^2(z;\psi^+)+(2\lambda-z)(2\mu-z)B(z;\psi^+)+z-2\lambda=0.
$$
Solving
\begin{equation*}\label{StAp}
B(z;\psi^+)=\frac{1}{2\lambda\mu}\left(2\mu-z-\sqrt{-\frac{z(2\mu-z)(2\lambda+2\mu-z)}{2\lambda-z}}\right).
\end{equation*}
If $\lambda\neq\mu$ the expression inside the square root is negative only for 
\begin{equation}\label{J1J2}
z\in J_1=[0,2\lambda\wedge2\mu]\quad\mbox{or}\quad z\in J_2=[2\lambda\vee2\mu,2\lambda+2\mu],
\end{equation}
where, as usual, $a\wedge b=\min\{a,b\}$ and $a\vee b=\max\{a,b\}$. It is also possible to see that there are no jumps. Therefore the spectral measure is given by
$$
\psi^+(x)=\frac{1}{2\pi\lambda\mu}\sqrt{\frac{x(2\mu-x)(2\lambda+2\mu-x)}{2\lambda-x}},\quad x\in J_1\cup J_2.
$$
The spectral measure for $\psi^-$ is the same but replacing $\lambda$ by $\mu$. Observe that if $\lambda=\mu$ we go back to the spectral measure of the absorbing $M/M/1$ queue in \eqref{mm1sp}.

Now, going back to the bilateral birth-death process and using the algebraic properties of $B(z;\psi^{\pm})$, we have, following \eqref{3BzZ} and after rationalizing, that
\begin{align*}
B(z;\psi_{11})&=-\sqrt{-\frac{2\mu-z}{z(2\lambda-z)(2\lambda+2\mu-z)}},\\
B(z;\psi_{12})&=-\frac{1}{2\lambda}\left(1+\sqrt{-\frac{(2\lambda-z)(2\mu-z)}{z(2\lambda+2\mu-z)}}\right),\\
B(z;\psi_{22})&=-\frac{\mu}{\lambda}\sqrt{-\frac{2\lambda-z}{z(2\mu-z)(2\lambda+2\mu-z)}}.
\end{align*}
Again, it is possible to see that there are no jumps. Therefore the spectral matrix has only an absolutely continuous part given by
\begin{equation*}
\Psi(x)=\frac{1}{\pi}\begin{pmatrix} \D\sqrt{\frac{2\mu-x}{x(2\lambda-x)(2\lambda+2\mu-x)}}&\D\frac{1}{2\lambda}\sqrt{\frac{(2\lambda-x)(2\mu-x)}{x(2\lambda+2\mu-x)}}\\ \D\frac{1}{2\lambda}\sqrt{\frac{(2\lambda-x)(2\mu-x)}{x(2\lambda+2\mu-x)}}&\D\frac{\mu}{\lambda}\sqrt{\frac{2\lambda-x}{x(2\mu-x)(2\lambda+2\mu-x)}}\end{pmatrix},\quad x\in J_1\cup J_2.
\end{equation*}
The explicit expression of the polynomials is now more difficult to compute. But using the main theorem in \cite{BGL} it is possible to obtain an explicit expression of the polynomials generated by the three-term recurrence relation \eqref{3TTRRZ}. If we define the new variable
$$
y=-1+\frac{(2\lambda-x)(2\mu-x)}{2\lambda\mu},
$$
then we have, for the first family,
\begin{align*}
Q_{2k}^1(x)&=(2y+1)U_{k-1}(y)-U_{k-2}(y),\quad k\in\mathbb{N},\quad Q_{2k+1}^1(x)=-\frac{1}{\lambda}(x-2\lambda)U_{k}(y),\quad k\in\mathbb{N}_0,\\
Q_{-2k-2}^1(x)&=-(2y+1)U_{k-1}(y)+U_{k-2}(y),\quad k\in\mathbb{N}_0,\quad Q_{-2k-1}^1(x)=\frac{1}{\lambda}(x-2\lambda)U_{k-1}(y),\quad k\in\mathbb{N},
\end{align*}
and for the second family we have $Q_1^2(x)=-1$ and
\begin{align*}
Q_{2k}^2(x)&=\frac{1}{\mu}(x-2\mu)U_{k-1}(y),\quad k\in\mathbb{N}_0,\quad Q_{2k+1}^2(x)=-(2y+1)U_{k-1}(y)+U_{k-2}(y),\quad k\in\mathbb{N},\\
Q_{-2k}^2(x)&=-\frac{1}{\mu}(x-2\mu)U_{k-1}(y),\quad k\in\mathbb{N}_0,\quad Q_{-2k-1}^2(x)=(2y+1)U_{k-1}(y)-U_{k-2}(y),\quad k\in\mathbb{N},
\end{align*}
where again $(U_n)_{n\in\mathbb{N}_0}$ are the Chebychev polynomials of the second kind. The transition probability functions $P_{ij}(t)$ can then be approximated using the Karlin-McGregor formula \eqref{3FKMcGZ}. 

\begin{figure*}[t!]
    \centering
    \begin{subfigure}[b]{0.5\textwidth}
        \centering
        \includegraphics[height=2.5in]{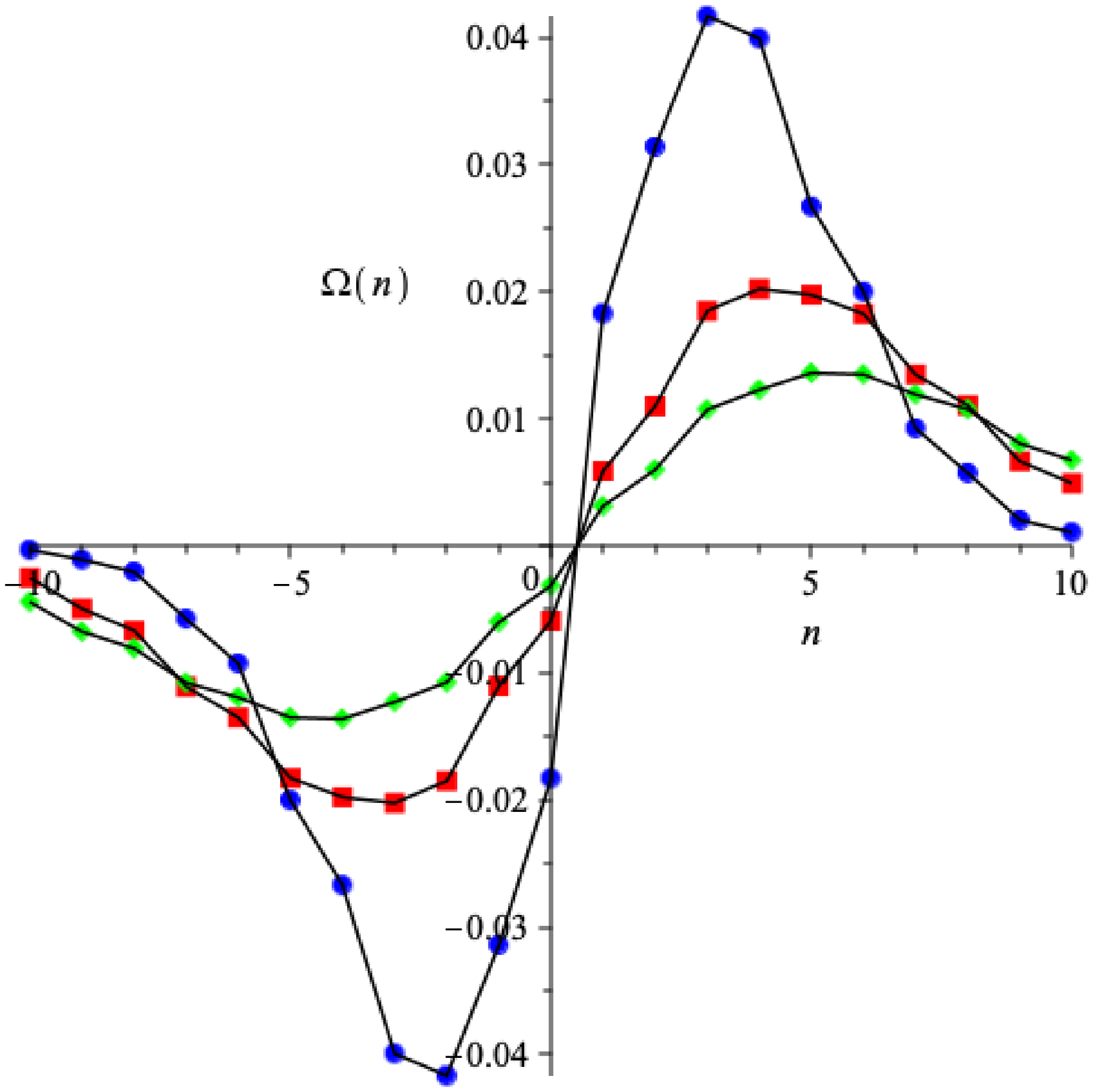}
        \caption{$\lambda=1,\mu=2$}
    \end{subfigure}%
    ~ 
    \begin{subfigure}[b]{0.5\textwidth}
        \centering
        \includegraphics[height=2.5in]{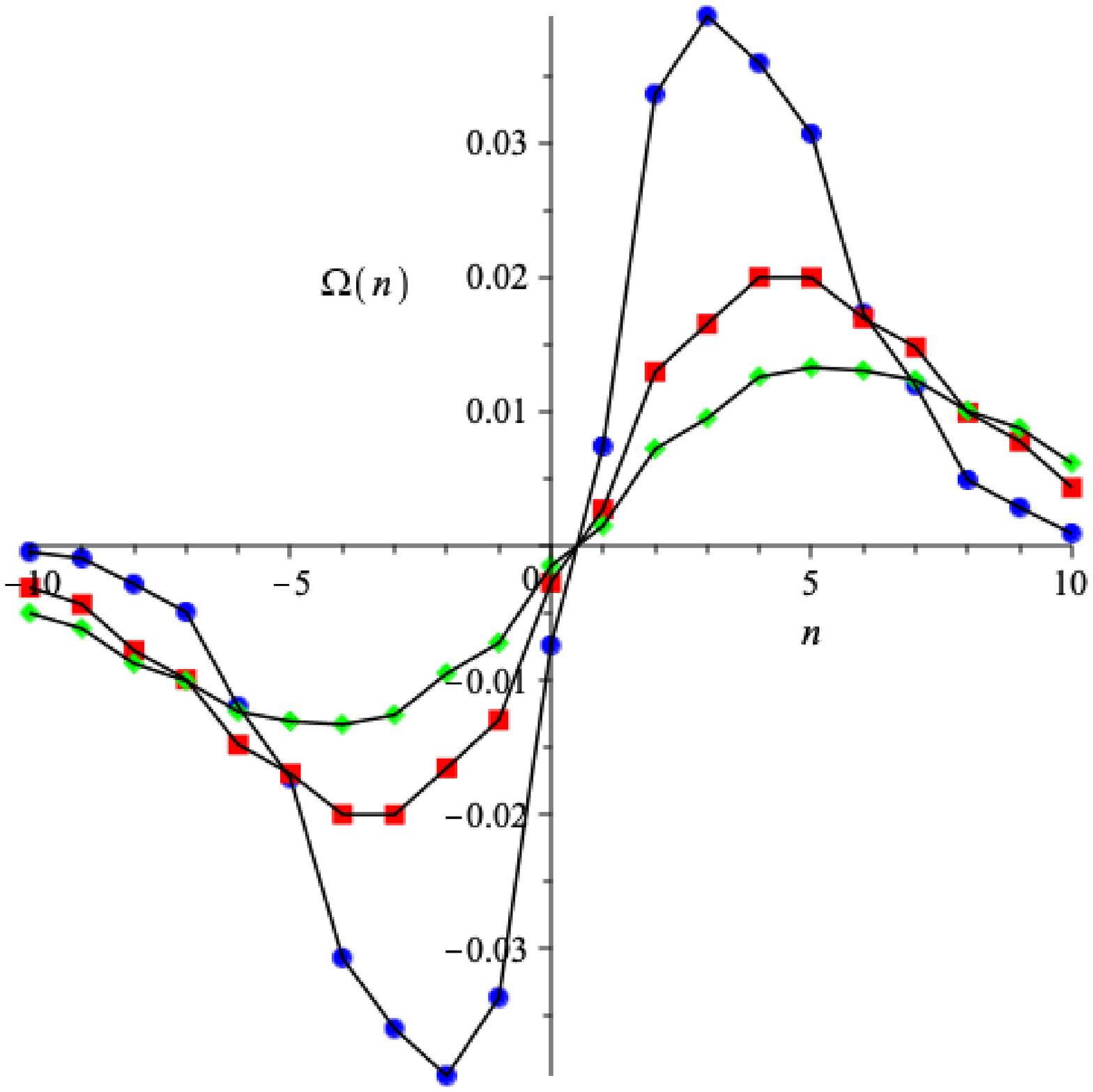}
        \caption{$\lambda=2,\mu=1$}
    \end{subfigure}
  \caption{The probability current \eqref{probcurZ} (here denoted by $\Omega(n)$) is plotted as function of $n$ for $j = 0$, $t = 3$ (blue circles), $t = 6$ (red squares) and $t =9$ (green diamonds).}  \label{fig3}
\end{figure*}

From \eqref{rec} and the explicit expression of the spectral matrix we have that $\int_0^\infty x^{-1}\psi_{11}(x)dx=\int_0^\infty x^{-1}\psi_{22}(x)dx=\infty$ for any values of $\lambda,\mu$. Therefore the process is always recurrent, as expected. Since there is no jump at the point 0, the process is always null recurrent. Again, we can get an approximation of the probability current by using \eqref{probcurZ}, the spectral matrix and the corresponding orthogonal polynomials. In Figure \ref{fig3} this probability current is plotted as a function of $n$ starting at $j=0$ for $t=3,6,9$ and for different values of the birth-death rates $\lambda,\mu$. 

\subsection{Case 2}

Consider the bilateral birth-death process with birth-death rates given by
$$
\lambda_{2n}=\lambda,\quad\lambda_{2n+1}=\mu, \quad \mu_{2n}=\mu,\quad \mu_{2n+1}=\lambda, \quad n\in\mathbb{Z},\quad \lambda,\mu>0. 
$$
The matrices $\mathcal{A}^{\pm}$ in \eqref{AAp} and \eqref{AAm} are now given by
\begin{equation*}
\mathcal{A}^+=\mathcal{A}^-=\begin{pmatrix}
-(\lambda+\mu)&\lambda&&&\\
\lambda&-(\lambda+\mu)&\mu&&\\
&\mu&-(\lambda+\mu)&\lambda&\\
&&\ddots&\ddots&\ddots
\end{pmatrix}.
\end{equation*}
Applying twice the identity (2.5) of \cite{KMc4} we obtain the following expression for the the Stieltjes transform of the measure $\psi^+$:
\begin{equation*}\label{StAp2}
B(z;\psi^+)=\frac{z^2-2(\lambda+\mu)z+2\mu(\lambda+\mu)-\sqrt{-z(2\mu-z)(2\lambda-z)(2\lambda+2\mu-z)}}{2\mu^2(\lambda+\mu-z)}.
\end{equation*}
Now $\psi^+$ will have an absolutely continuous part and a discrete part. Indeed
$$
\psi^+(x)=\frac{\sqrt{x(2\lambda-x)(2\mu-x)(2\lambda+2\mu-x)}}{2\pi\mu^2|\lambda+\mu-x|}\mathbf{1}_{\{J_1\cup J_2\}},+\left(1-\frac{\lambda^2}{\mu^2}\right)\delta_{\lambda+\mu}(x)\mathbf{1}_{\{\mu>\lambda\}},
$$
where $\mathbf{1}_A$ is the indicator function, $\delta_a(x)$ is the Dirac delta located at $x=a$ and $J_1,J_2$ are defined by \eqref{J1J2}. Now, going back to the bilateral birth-death process and using the algebraic properties of $B(z;\psi^{+})=B(z;\psi^{-})$, we have, following \eqref{3BzZ} and after rationalizing, that
\begin{align*}
B(z;\psi_{11})&=B(z;\psi_{22})=-\frac{\lambda+\mu-z}{\sqrt{-z(2\mu-z)(2\lambda-z)(2\lambda+2\mu-z)}},\\
B(z;\psi_{12})&=-\frac{1}{2\mu}\left(1+\frac{(\lambda+\mu-z)^2+\mu^2-\lambda^2}{\sqrt{-z(2\mu-z)(2\lambda-z)(2\lambda+2\mu-z)}}\right),
\end{align*}
Again, it is possible to see that there are no jumps. Therefore the spectral matrix has only an absolutely continuous part given by
\begin{equation*}
\Psi(x)=\frac{1}{\pi\sqrt{x(2\mu-x)(2\lambda-x)(2\lambda+2\mu-x)}}\begin{pmatrix} |\lambda+\mu-x|&|(\lambda+\mu-x)^2+\mu^2-\lambda^2|\\ |(\lambda+\mu-x)^2+\mu^2-\lambda^2|&|\lambda+\mu-x|\end{pmatrix},
\end{equation*}
where $x\in J_1\cup J_2$ and $J_1,J_2$ are defined by \eqref{J1J2}. Again, by using the main theorem in \cite{BGL}, it is possible to obtain an explicit expression of the polynomials generated by the three-term recurrence relation \eqref{3TTRRZ}. If we define the new variable
$$
y=\frac{(\lambda+\mu-x)^2-\lambda^2-\mu^2}{2\lambda\mu},
$$
then we have $Q_1^2(x)=-\mu/\lambda$ and
\begin{align*}
Q_{2k}^1(x)&=(2y+\mu/\lambda)U_{k-1}(y)-U_{k-2}(y),\quad k\in\mathbb{N},\quad Q_{2k+1}^1(x)=\frac{1}{\lambda}(\lambda+\mu-x)U_{k}(y),\quad k\in\mathbb{N}_0,\\
Q_{2k}^2(x)&=\frac{1}{\lambda}(x-\lambda-\mu)U_{k-1}(y),\quad k\in\mathbb{N}_0,\quad Q_{2k+1}^2(x)=-\frac{\mu}{\lambda}\left[(2y+\lambda/\mu)U_{k-1}(y)-U_{k-2}(y)\right],\quad k\in\mathbb{N},
\end{align*}
while for the negative indices we have $Q_{-n-1}^1(x)=Q_n^2(x)$ and $Q_{-n-1}^2(x)=Q_n^1(x)$ for $n\in\mathbb{N}_0$. Again $(U_n)_{n\in\mathbb{N}_0}$ are the Chebychev polynomials of the second kind. The transition probability functions $P_{ij}(t)$ can then be approximated using the Karlin-McGregor formula \eqref{3FKMcGZ}. 
The probabilistic properties of this process are the same as in the previous case, i.e. the process is always null recurrent. Also graphs for the probability current are similar as the ones plotted in Figure \ref{fig3}.

\section{Variants of the bilateral birth-death processes with constant rates}\label{sec4}

In this section we will study a couple of variants of the bilateral birth-death processes studied in Examples \ref{exb1} and \ref{exb2}, allowing one defect at the state 0. Although we are introducing a small change we will see that the computations get more involved than usual.

\subsection{Case 1}

The birth-death rates are now given by
$$
\lambda_n=\lambda,\quad\mu_n=\mu,\quad n\in\mathbb{Z}\setminus\{0\},\quad\lambda,\mu,\lambda_0,\mu_0>0.
$$
Now the matrix $\mathcal{A}^+$ in \eqref{AAp} is given by 
\begin{equation*}\label{varA1}
\mathcal{A}^+=\begin{pmatrix}
-(\lambda_0+\mu_0)&\lambda_0&&&\\
\mu&-(\lambda+\mu)&\lambda&&\\
&\mu&-(\lambda+\mu)&\lambda&\\
&&\ddots&\ddots&\ddots
\end{pmatrix},
\end{equation*}
while the matrix $\mathcal{A}^-$ in \eqref{AAm} is the same as the one of the absorbing $M/M/1$ queue in Example \ref{3Ej7} and also the infinitesimal operator of the 0-th birth-death process associated to $\mathcal{A}^+$ (i.e. the infinitesimal operator defined from $\mathcal{A^+}$ eliminating the first row and column of $\mathcal{A^+}$). Calling $B(z)$ the Stieltjes transform in \eqref{31server1} and Applying (2.5) of \cite{KMc4} we obtain
\begin{equation}\label{StAp3}
B(z;\psi^+)=\frac{1}{\lambda_0+\mu_0-z-\lambda_0\mu B(z)},\quad B(z;\psi^-)=B(z).
\end{equation}
Following \eqref{3BzZ} and rationalizing, we obtain that
\begin{equation*}
B(z;\psi_{\alpha\beta})=\frac{p_{\alpha\beta}(z)+q_{\alpha\beta}(z)\sqrt{(\lambda+\mu-z)^2-4\lambda\mu}}{D(z)},\quad \alpha,\beta=1,2,
\end{equation*}
where 
\begin{align*}
p_{11}(z)&=(\lambda_0\mu+\mu_0\lambda-2\lambda\mu)z+(\lambda-\mu)(\lambda_0\mu-\mu_0\lambda),\quad q_{11}(z)=-(\lambda_0\mu+\lambda\mu_0),\\
p_{12}(z)&=\lambda\left[z^2-(\lambda_0+\mu_0+\lambda+\mu)z+(\lambda - \mu)(\lambda_0 - \mu_0)\right],\quad q_{12}(z)=-\lambda(\lambda_0+\mu_0-z),\\
p_{22}(z)&=\frac{1}{\mu_0}\left[(\lambda_0 - \lambda)z^3+(-\lambda_0^2 - \lambda_0\mu_0 - 2\lambda_0\mu + 2\mu_0\lambda + \lambda^2 + \lambda\mu)z^2+\mu_0(\lambda - \mu)(\lambda_0\mu - \mu_0\lambda)\right.\\
&\qquad \left.+(\lambda_0^2\lambda + \lambda_0^2\mu - \lambda_0\mu_0\lambda + 2\lambda_0\mu_0\mu - \lambda_0\lambda^2 + \lambda_0\mu^2 - \mu_0^2\lambda - 2\mu_0\lambda\mu)z\right],\\
q_{22}(z)&=-\frac{1}{\mu_0}\left[(\lambda-\lambda_0)z^2+(\lambda_0^2 + \lambda_0\mu_0 - \lambda_0\lambda + \lambda_0\mu - 2\mu_0\lambda)z+\mu_0(2\lambda_0\lambda - \lambda_0\mu + \mu_0\lambda)\right],\\
D(z)&=(-2\lambda_0\mu-2\mu_0\lambda+2\lambda\mu)z^2-2\lambda_0\mu_0(\lambda - \mu)^2\\
&\qquad +(2\lambda_0^2\mu + 2\lambda_0\mu_0\lambda + 2\lambda_0\mu_0\mu - 2\lambda_0\lambda\mu + 2\lambda_0\mu^2 + 2\mu_0^2\lambda + 2\mu_0\lambda^2 - 2\mu_0\lambda\mu)z.
\end{align*}
The spectral matrix $\Psi(x)=\Psi_c(x)+\Psi_d(x)$ has now an absolutely continuous part $\Psi_c(x)$, given by 
\begin{equation*}
\Psi(x)=\frac{\sqrt{(x-\sigma_-)(\sigma_+-x)}}{\pi D(x)}
\begin{pmatrix}\lambda_0\mu+\lambda\mu_0&\lambda(\lambda_0+\mu_0-x)\\\lambda(\lambda_0+\mu_0-x)&q_{22}(x)\end{pmatrix},\quad x\in[\sigma_-,\sigma_+],
\end{equation*}
where $\sigma_{\pm}=(\sqrt{\lambda}\pm\sqrt{\mu})^2$. For the discrete part $\Psi_d(x)$ we need to study the poles of $B(z;\psi_{\alpha\beta}), \alpha,\beta=1,2,$ which are the roots of the second-degree polynomial $D(z)$. For that let us introduce some constants which will simplify considerably the sequel:
\begin{align*}
R&=(\lambda_0 - \mu_0 - \lambda + \mu)^2 + 4\lambda_0\mu_0,\\
C&=\lambda\mu R-(\lambda_0\mu+\mu_0\lambda-\lambda\mu)(\lambda_0\mu+\mu_0\lambda+2\lambda\mu).
\end{align*}
The roots of $D(z)$ are then given by
\begin{equation}\label{gammas}
\gamma_{\pm}=\frac{(\lambda_0\mu + \mu_0\lambda)(\lambda_0 + \mu_0) + (\lambda_0\mu - \mu_0\lambda)(\mu-\lambda)\pm (\lambda_0\mu + \mu_0\lambda)\sqrt{R}}{2(\lambda_0\mu + \mu_0\lambda - \lambda\mu)}.
\end{equation}
Observe that $\sqrt{R}$ is well-defined since $R>0$. A long but straightforward computation gives the magnitudes for each of these poles. Defining the constants
\begin{align*}
A_{11}^{\pm}&=\frac{1}{2(\lambda_0\mu+\mu_0\lambda-\lambda\mu)}\left[\lambda_0\mu+\mu_0\lambda-2\lambda\mu\pm(\lambda_0\mu+\mu_0\lambda)(\lambda+\mu-\lambda_0-\mu_0)\frac{\sqrt{R}}{R}\right],\\
A_{12}^{\pm}&=\frac{\lambda}{2(\lambda_0\mu+\mu_0\lambda-\lambda\mu)^2}\left[\lambda\mu(\lambda+\mu-\lambda_0-\mu_0)\pm(C-(\lambda_0\mu+\mu_0\lambda-\lambda\mu)(\lambda_0\mu+\mu_0\lambda-2\lambda\mu))\frac{\sqrt{R}}{R}\right],\\
A_{22}^{\pm}&=\frac{\lambda^2}{2(\lambda_0\mu+\mu_0\lambda-\lambda\mu)^3}\left[-C\pm(\lambda+\mu-\lambda_0-\mu_0)(C+2\lambda\mu(\lambda_0\mu+\mu_0\lambda-\lambda\mu))\frac{\sqrt{R}}{R}\right],
\end{align*}
we have that the discrete part $\Psi_d(x)$ is given by
\begin{equation*}
\Psi_d(x)=\begin{pmatrix}A_{11}^-&A_{12}^-\\A_{12}^-&A_{22}^-\end{pmatrix}\delta_{\gamma_+}(x)\mathbf{1}_{\{A_{11}^->0,\gamma_+>0\}}+\begin{pmatrix}A_{11}^+&A_{12}^+\\A_{12}^+&A_{22}^+\end{pmatrix}\delta_{\gamma_-}(x)\mathbf{1}_{\{A_{11}^+>0,\gamma_->0\}},
\end{equation*}
where $\mathbf{1}_A$ is the indicator function and $\delta_a(x)$ is the Dirac delta located at $x=a$. We have not been able to find a simplification of the conditions $\{A_{11}^->0,\gamma_+>0\}$ and $\{A_{11}^+>0,\gamma_->0\}$ in terms of the birth-death rates $\lambda,\mu,\lambda_0,\mu_0$. Finally, we can also compute the polynomials generated by the three-term recurrence relation \eqref{3TTRRZ}. If we define the new variable
\begin{equation}\label{yyx1}
y=\frac{\lambda+\mu-x}{2\sqrt{\lambda\mu}},
\end{equation}
then we have
\begin{align}
\nonumber Q_n^1(x)&=\left(\frac{\mu}{\lambda}\right)^{n/2}\left[\frac{2(\lambda_0-\lambda)}{\lambda_0}T_n\left(y\right)+\frac{2\lambda-\lambda_0}{\lambda_0}U_n\left(y\right)+\sqrt{\frac{\lambda}{\mu}}\frac{\lambda_0+\mu_0-\lambda-\mu}{\lambda_0}U_{n-1}\left(y\right)\right],\quad n\in\mathbb{N}_0,\\
\label{QQQ1}Q_{-n-1}^1(x)&=-\left(\frac{\lambda}{\mu}\right)^{(n+1)/2}U_{n-1}\left(y\right),\quad n\in\mathbb{N}_0,\\
\nonumber Q_n^2(x)&=-\frac{\mu_0}{\lambda_0}\left(\frac{\mu}{\lambda}\right)^{(n-1)/2}U_{n-1}\left(y\right),\quad Q_{-n-1}^2(x)=\left(\frac{\lambda}{\mu}\right)^{n/2}U_{n}\left(y\right),\quad n\in\mathbb{N}_0,
\end{align}
where $(T_n)_{n\in\mathbb{N}_0}$ and $(U_n)_{n\in\mathbb{N}_0}$ are the Chebychev polynomials of the first and second kind, respectively. The transition probability functions $P_{ij}(t)$ can then be approximated using the Karlin-McGregor formula \eqref{3FKMcGZ}. It is possible to see that if $\lambda_0=\lambda$ and $\mu_0=\mu$ we go back to the Example \ref{exb1} of the bilateral birth-death process with constant rates.

\begin{figure*}[t!]
    \centering
    \begin{subfigure}[b]{0.5\textwidth}
        \centering
        \includegraphics[height=2.5in]{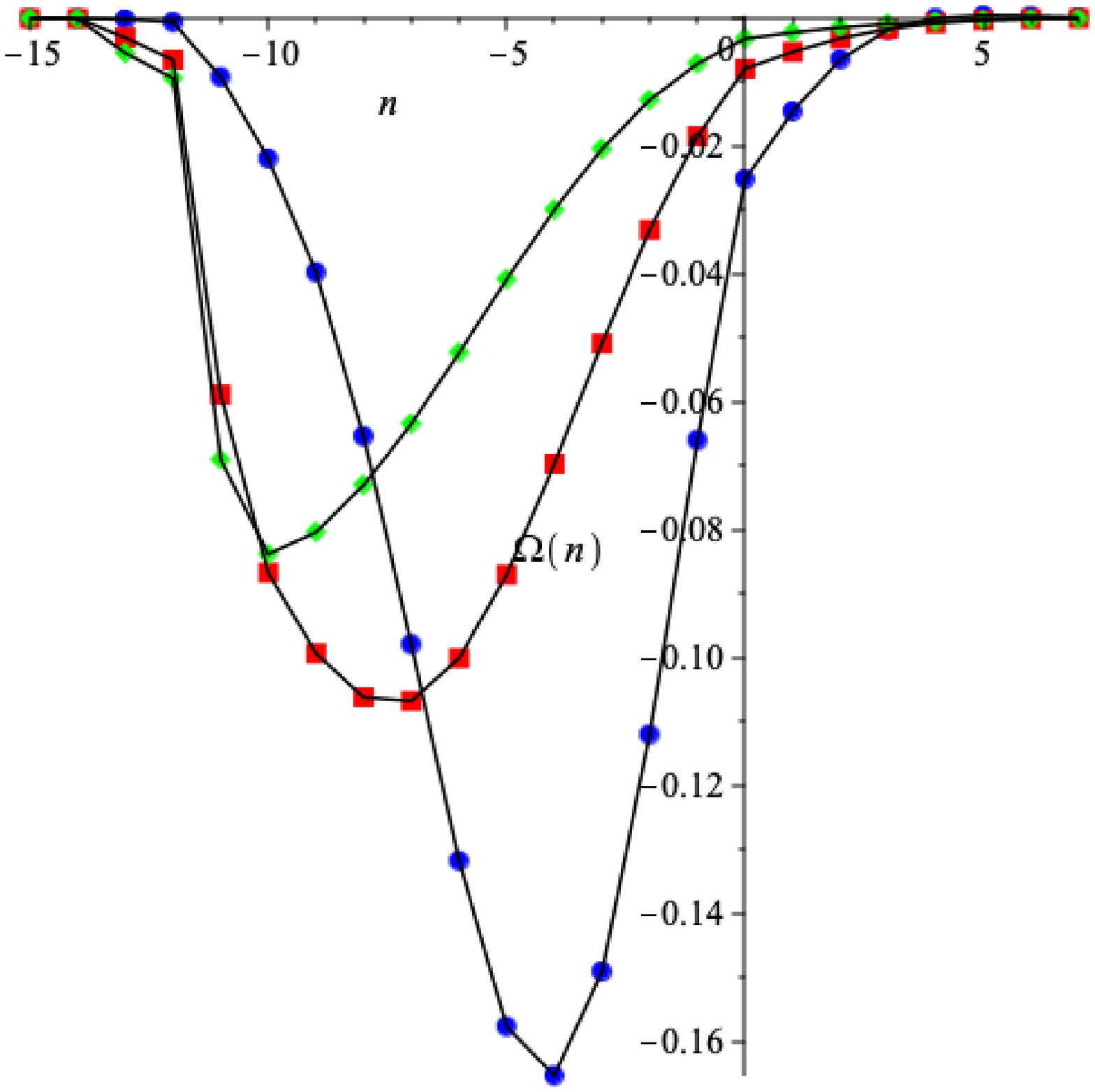}
        \caption{$\lambda=1,\mu=2,\lambda_0=1,\mu_0=5$}
    \end{subfigure}%
    ~ 
    \begin{subfigure}[b]{0.5\textwidth}
        \centering
        \includegraphics[height=2.5in]{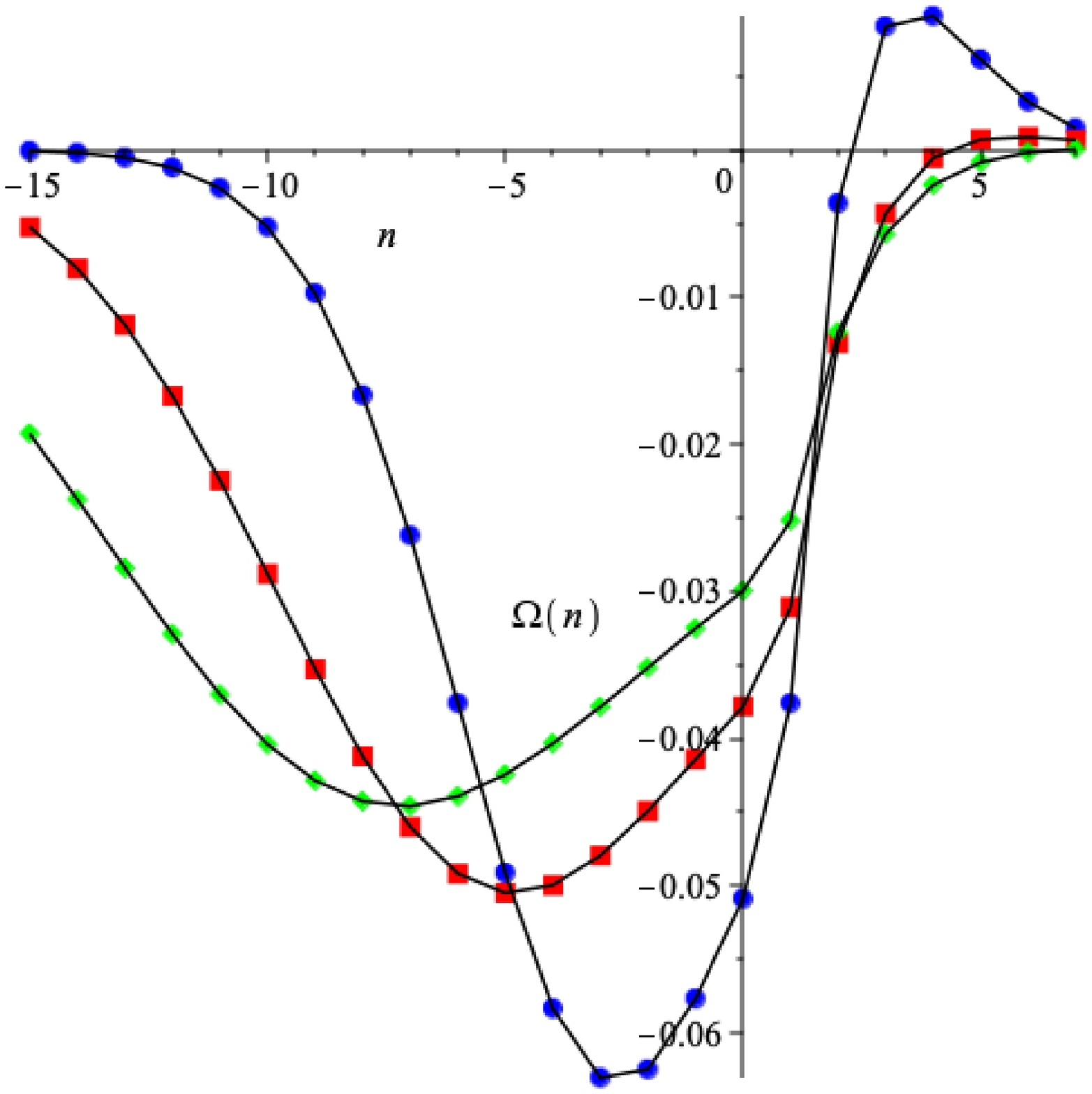}
        \caption{$\lambda=1,\mu=2,\lambda_0=5,\mu_0=1$}
    \end{subfigure}
    \\
      \centering
    \begin{subfigure}[b]{0.5\textwidth}
        \centering
        \includegraphics[height=2.5in]{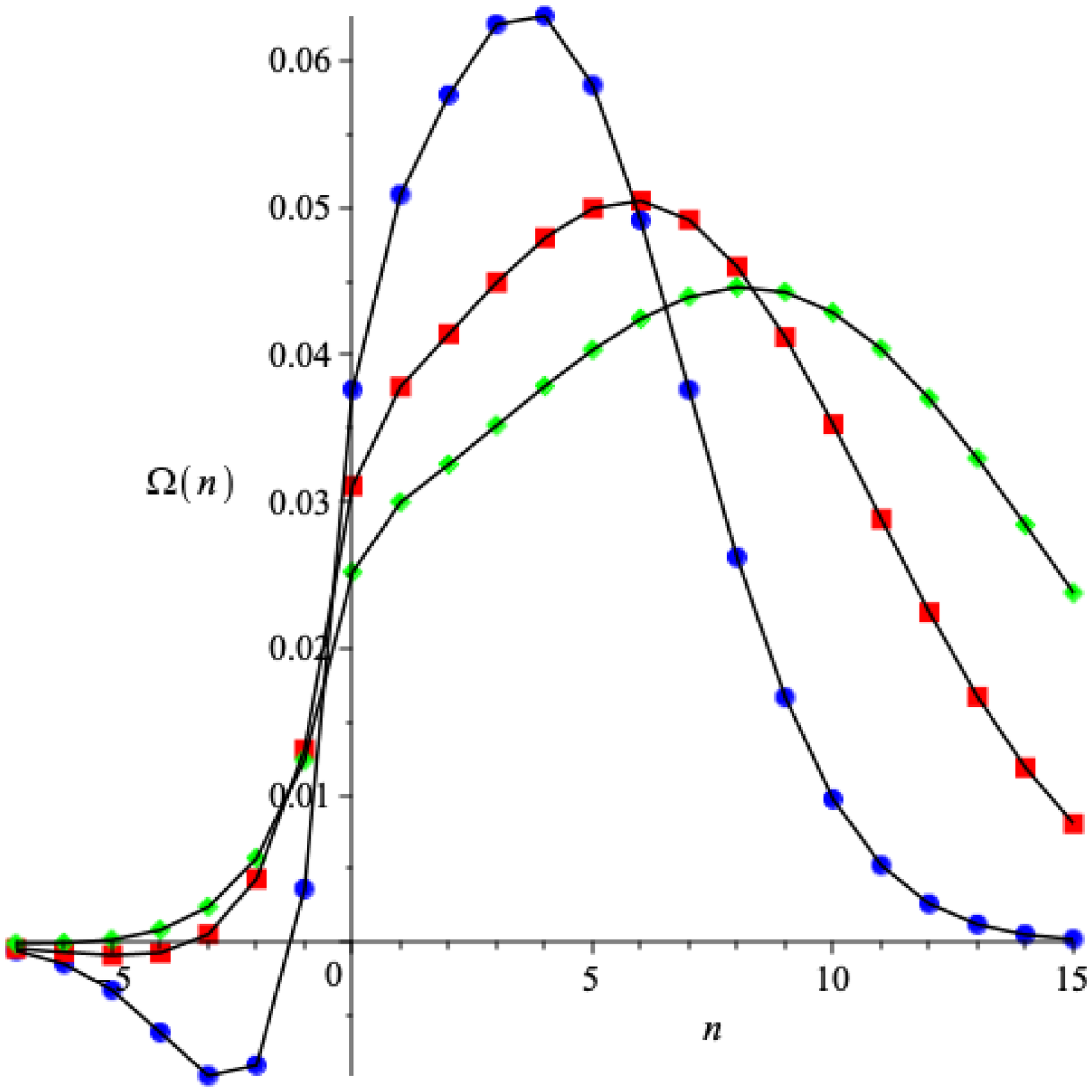}
        \caption{$\lambda=2,\mu=1,\lambda_0=1,\mu_0=5$}
    \end{subfigure}%
    ~ 
    \begin{subfigure}[b]{0.5\textwidth}
        \centering
        \includegraphics[height=2.5in]{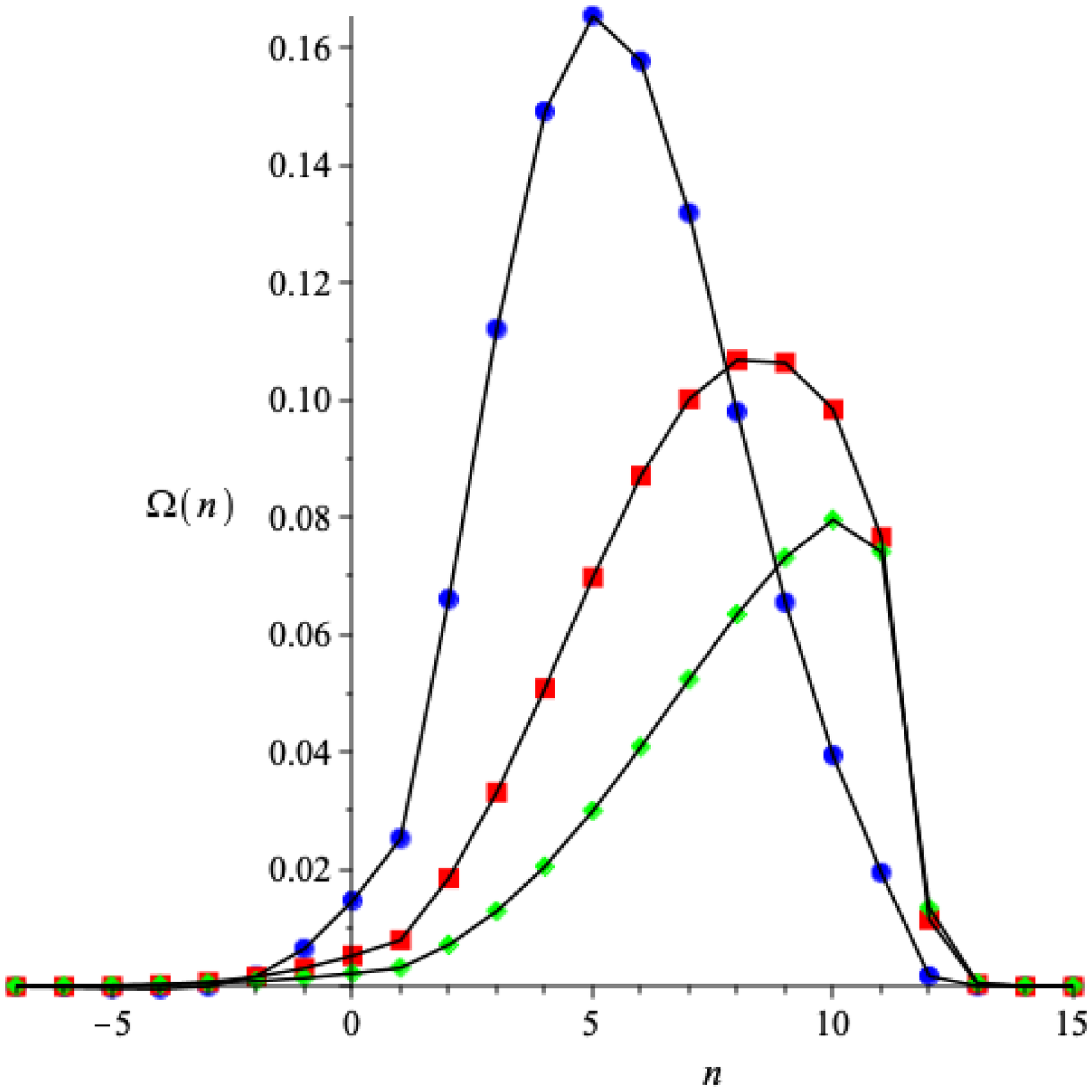}
        \caption{$\lambda=2,\mu=1,\lambda_0=5,\mu_0=1$}
    \end{subfigure}
      \caption{The probability current \eqref{probcurZ} (here denoted by $\Omega(n)$) is plotted as function of $n$ for $j = 0$, $t = 3$ (blue circles), $t = 6$ (red squares) and $t =9$ (green diamonds).}  \label{fig4}
\end{figure*}

From \eqref{rec} and the explicit expression of the spectral matrix we have that $\int_0^\infty x^{-1}\psi_{11}(x)dx<\infty$ and $\int_0^\infty x^{-1}\psi_{22}(x)dx<\infty$ for $\lambda\neq\mu$ and any values of $\lambda_0,\mu_0$. Therefore the process is always transient unless $\lambda=\mu$, where it is recurrent. In that case we have that $\gamma_{\pm}$ in \eqref{gammas} are given by $\gamma_-=0$ and $\gamma_+=(\lambda_0+\mu_0)^2/(\alpha+\beta-\lambda)$ and $A_{ij}^-=0,i,j=1,2$. Therefore there is no jump at the point 0 and the process is null recurrent. Again, we can get an approximation of the probability current by using \eqref{probcurZ}, the spectral matrix and the corresponding orthogonal polynomials. In Figure \ref{fig4} this probability current is plotted as a function of $n$ starting at $j=0$ for $t=3,6,9$ and for different values of the birth-death rates $\lambda,\mu,\lambda_0,\mu_0$. In the first and fourth cases we only have one Dirac delta located at $x=\gamma_+$, while in the the second and third cases we have two Dirac deltas located at $x=\gamma_{\pm}$.

\subsection{Case 2}

Let us consider now a variant of the symmetric bilateral birth-death process with constant rates studied in Example \ref{exb2}. The birth-death rates are now given by
$$
\lambda_n=\lambda,\quad\mu_n=\mu,\quad \lambda_{-n}=\mu,\quad \mu_{-n}=\lambda, \quad n\in\mathbb{N},\quad \lambda,\mu,\lambda_0,\mu_0>0. 
$$
The matrices $\mathcal{A}^{\pm}$ in \eqref{AAp} and \eqref{AAm} are the same as in the previous case, so the we can use the same notation for the Stieltjes transforms of $\psi^{\pm}$ in \eqref{StAp3}. Now there is a small change in formulas \eqref{3BzZ} but this somehow simplify a little bit the computation of the Stieltjes transforms $B(z;\psi_{\alpha\beta}), \alpha,\beta=1,2$. Again, after rationalizing, we obtain that
\begin{equation*}
B(z;\psi_{\alpha\beta})=\frac{p_{\alpha\beta}(z)+q_{\alpha\beta}(z)\sqrt{(\lambda+\mu-z)^2-4\lambda\mu}}{D(z)},\quad \alpha,\beta=1,2,
\end{equation*}
where 
\begin{align*}
p_{11}(z)&=(\lambda_0+\mu_0-2\lambda)z+(\lambda -\mu)(\lambda_0 +\mu_0),\quad q_{11}(z)=-(\lambda_0+\mu_0),\\
p_{12}(z)&=z^2-(\lambda_0 + \mu_0 + \lambda + \mu)z+(\lambda - \mu)(\lambda_0 + \mu_0),\quad q_{12}(z)=-(\lambda_0+\mu_0-z),\\
p_{22}(z)&=\frac{1}{\lambda\mu_0}\left[(\lambda_0 - \lambda)z^3+(-\lambda_0^2 - \lambda_0\mu_0 - 2\lambda_0\mu + 2\mu_0\lambda + \lambda^2 + \lambda\mu)z^2+\mu_0\lambda(\lambda - \mu)(\lambda_0 + \mu_0)\right.\\
&\qquad \left.+(\lambda_0^2\lambda + \lambda_0^2\mu + \lambda_0\mu_0\mu - \lambda_0\lambda^2 + \lambda_0\mu^2 - \mu_0^2\lambda - 2\mu_0\lambda^2)z\right],\\
q_{22}(z)&=-\frac{1}{\lambda\mu_0}\left[(\lambda-\lambda_0)z^2+(\lambda_0^2 + \lambda_0\mu_0 - \lambda_0\lambda + \lambda_0\mu - 2\mu_0\lambda)z+\mu_0\lambda(\lambda_0 + \mu_0)\right],\\
D(z)&=2z[(\lambda_0 +\mu_0)(\lambda_0 + \mu_0 - \lambda + \mu)-(\lambda_0+\mu_0-\lambda)z].
\end{align*}
The spectral matrix $\Psi(x)=\Psi_c(x)+\Psi_d(x)$ has now an absolutely continuous part $\Psi_c(x)$, given by 
\begin{equation}\label{wvar1}
\Psi_c(x)=\frac{\sqrt{(x-\sigma_-)(\sigma_+-x)}}{\pi D(x)}
\begin{pmatrix}\lambda_0+\mu_0&\lambda_0+\mu_0-x\\\lambda_0+\mu_0-x&q_{22}(x)\end{pmatrix},\quad x\in[\sigma_-,\sigma_+],
\end{equation}
where $\sigma_{\pm}=(\sqrt{\lambda}\pm\sqrt{\mu})^2$. For the discrete part $\Psi_d(x)$ we get now easier expressions since the roots of $D(z)$ are given by 0 and the constant
\begin{equation*}
\eta=\frac{(\lambda_0+\mu_0)(\lambda_0+\mu_0-\lambda+\mu)}{\lambda_0+\mu_0-\lambda}.
\end{equation*}
After some straightforward computations we have that the discrete part $\Psi_d(x)$ is given by
\begin{align}
\label{wvar2}\Psi_d(x)&=\frac{\mu-\lambda}{\lambda_0+\mu_0+\mu-\lambda}\begin{pmatrix}1&1\\1&1\end{pmatrix}\delta_{0}(x)\mathbf{1}_{\{\mu>\lambda\}}\\
\nonumber&+\frac{(\lambda_0+\mu_0-\lambda)^2-\lambda\mu}{(\lambda_0+\mu_0-\lambda)(\lambda_0+\mu_0+\mu-\lambda)}\begin{pmatrix}1&-\D\frac{\mu}{\lambda_0+\mu_0-\lambda}\\-\D\frac{\mu}{\lambda_0+\mu_0-\lambda}&\D\left(\frac{\mu}{\lambda_0+\mu_0-\lambda}\right)^2\end{pmatrix}\delta_{\eta}(x)\mathbf{1}_{\{|\lambda_0+\mu_0-\lambda|>\sqrt{\lambda\mu}\}},
\end{align}
where $\mathbf{1}_A$ is the indicator function and $\delta_a(x)$ is the Dirac delta located at $x=a$. Finally, we can also compute the polynomials generated by the three-term recurrence relation \eqref{3TTRRZ}. Using the same notation as in \eqref{yyx1} we have that $Q_n^1(x)$ and $Q_n^2(x)$ for $n\in\mathbb{N}_0$ are the same as in the previous case in \eqref{QQQ1}, while for $Q_{-n-1}^1(x)$ and $Q_{-n-1}^2(x)$ we only have to change $\lambda$ by $\mu$ in \eqref{QQQ1}. Once more, the transition probability functions $P_{ij}(t)$ can then be approximated using the Karlin-McGregor formula \eqref{3FKMcGZ}. It is possible to see that if $\lambda_0=\lambda$ and $\mu_0=\mu$ we go back to the Example \ref{exb2} of the symmetric bilateral birth-death process with constant rates.

\begin{figure*}[t!]
    \centering
    \begin{subfigure}[b]{0.5\textwidth}
        \centering
        \includegraphics[height=2.5in]{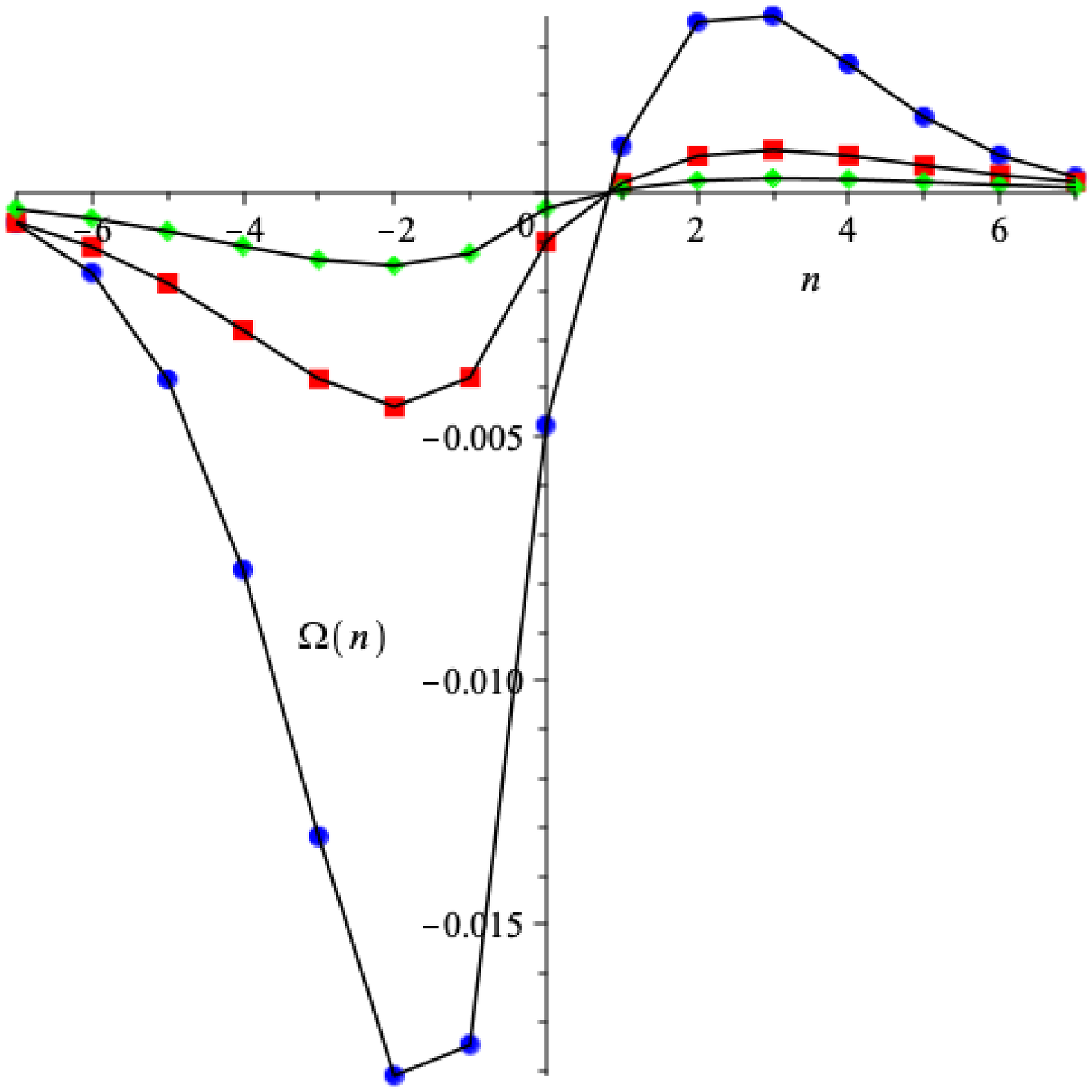}
        \caption{$\lambda=1,\mu=2,\lambda_0=1,\mu_0=5$}
    \end{subfigure}%
    ~ 
    \begin{subfigure}[b]{0.5\textwidth}
        \centering
        \includegraphics[height=2.5in]{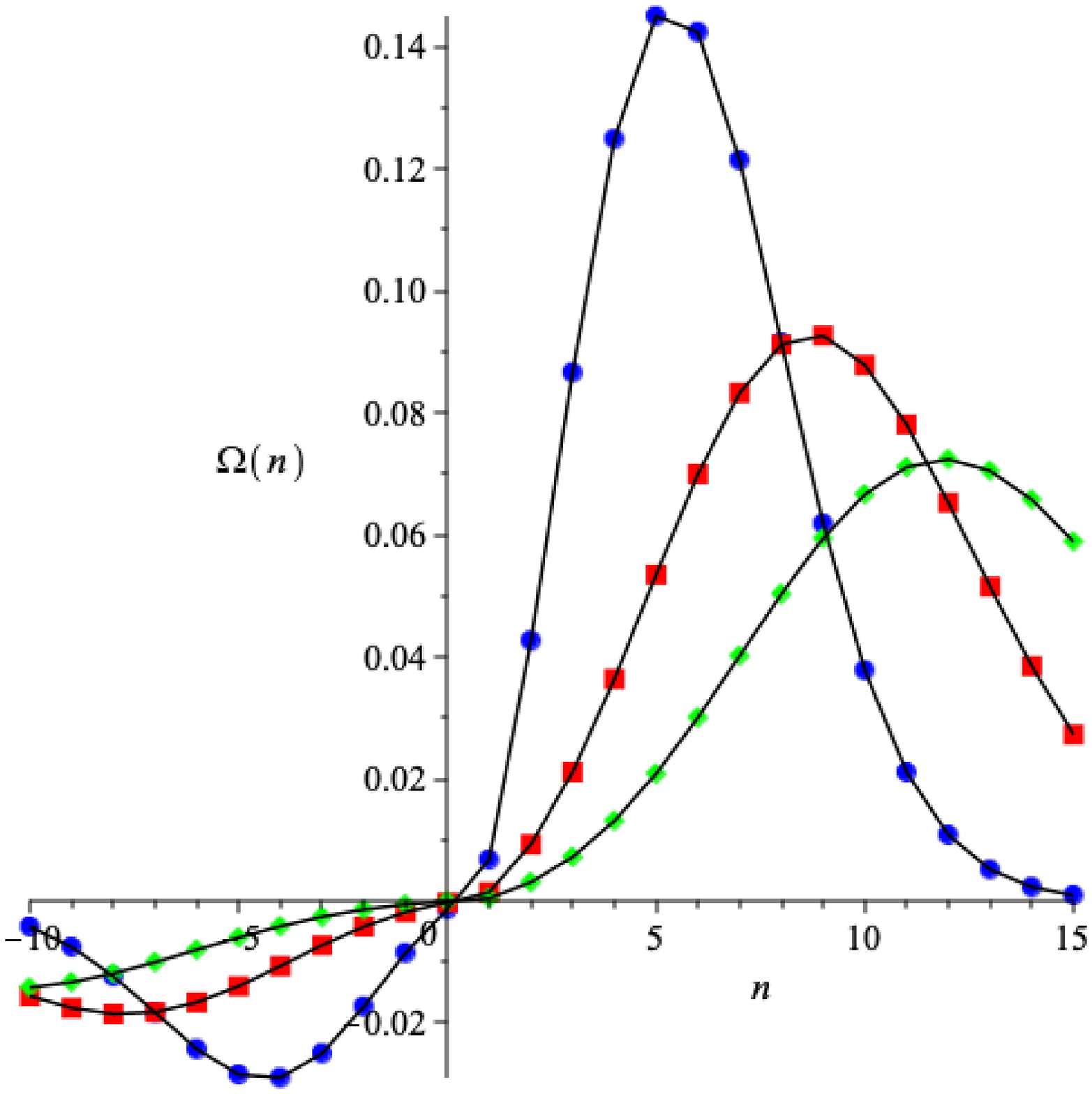}
        \caption{$\lambda=2,\mu=1,\lambda_0=5,\mu_0=1$}
    \end{subfigure}
  \caption{The probability current \eqref{probcurZ} (here denoted by $\Omega(n)$) is plotted as function of $n$ for $j = 0$, $t = 3$ (blue circles), $t = 6$ (red squares) and $t =9$ (green diamonds).}  \label{fig6}
\end{figure*}

From \eqref{rec} and the explicit expression of the spectral matrix we have that $\int_0^\infty x^{-1}\psi_{11}(x)dx<\infty$ and $\int_0^\infty x^{-1}\psi_{22}(x)dx<\infty$ for $\lambda>\mu$ and any values of $\lambda_0,\mu_0$. Therefore the process is always transient for $\lambda>\mu$. If $\lambda\leq\mu$ then both integrals diverge and the process is recurrent. For $\lambda<\mu$ we always have a jump at the point 0, so the process is positive recurrent. If $\lambda=\mu$ then the process is null recurrent. Since the potential coefficients are given here by (see \eqref{3potcoefZ})
$$
\pi_0=1,\quad\pi_n=\frac{\lambda_0}{\mu}\left(\frac{\lambda}{\mu}\right)^{n-1},\quad\pi_{-n}=\frac{\mu_0}{\mu}\left(\frac{\lambda}{\mu}\right)^{n-1},\quad n\in\mathbb{N},
$$
we have that the invariant distribution $\bm\pi$ for this process is given by
$$
\bm\pi=\frac{\mu-\lambda}{\mu-\lambda+\mu_0+\lambda_0}\left(\cdots,\frac{\mu_0\lambda}{\mu^2},\frac{\mu_0}{\mu},1,\frac{\lambda_0}{\mu},\frac{\lambda_0\lambda}{\mu^2},\frac{\lambda_0\lambda^2}{\mu^3},\cdots\right),\quad \mu>\lambda.
$$
A couple of graphs of the probability current \eqref{probcurZ} can be found in Figure \ref{fig6}. In the first case we have two Dirac deltas located at 0 and $42/5$ while in the second case we only have one located at $x=15/2$.

\begin{remark}
In Section 5.3 of \cite{GN} it was studied this example for the particular case of $\lambda_0=\mu_0=\lambda$. Substituting these values in \eqref{wvar1} and \eqref{wvar2} we have that the spectral matrix in this particular case is given by $\Psi(x)=\Psi_c(x)+\Psi_d(x)$ where
\begin{equation*}
\Psi_c(x)=\frac{\sqrt{(x-\sigma_-)(\sigma_+-x)}}{\pi x(2\lambda+2\mu-x)}\begin{pmatrix}1&\D 1-\frac{x}{2\lambda}\\
\D 1-\frac{x}{2\lambda}&-1+\D\left(1-\frac{\mu}{\lambda}\right)\frac{x}{2\lambda}\end{pmatrix},\quad x\in[\sigma_-,\sigma_+],
\end{equation*}
where $\sigma_{\pm}=(\sqrt{\lambda}\pm\sqrt{\mu})^2$ and
\begin{equation*}
\Psi_d(x)=\frac{\mu-\lambda}{\lambda+\mu}\left[\begin{pmatrix}1&1\\1&1\end{pmatrix}\delta_0(x)\mathbf{1}_{\{\mu>\lambda\}}+\begin{pmatrix}1&-\mu/\lambda\\-\mu/\lambda&(\mu/\lambda)^2\end{pmatrix}\delta_{2\lambda+2\mu}(x)\mathbf{1}_{\{\lambda>\mu\}}\right].
\end{equation*}
An extensive analysis of the probability properties of this example was given in Section 5.3 of \cite{GN}.
\end{remark}

\section{Splitting into two different $M/M/1$ queues}\label{sec5}

In this section we will study a bilateral birth-death process where $\mathcal{A}^{\pm}$ are two different absorbing $M/M/1$ queues with constant birth-death rates given by $\lambda,\mu>0$ and $\beta,\alpha>0$, respectively. Therefore we have
$$
\lambda_n=\lambda,\quad\mu_n=\mu,\quad \lambda_{-n-1}=\beta,\quad\mu_{-n-1}=\alpha,\quad n\in\mathbb{N}_0,\quad\lambda,\mu,\alpha,\beta>0.
$$
The Stieltjes transforms $B(z;\psi^{\pm})$ are given by \eqref{31server1} (replacing $\lambda,\mu$ by $\alpha,\beta$ in $B(z;\psi^{-})$). For simplicity we will use the following notation
\begin{equation}\label{Sigmab}
\Sigma^{a,b}_{\pm}(z)=\sqrt{\pm\left[(a+b-z)^2-4ab\right]},
\end{equation}
By using algebraic properties of $B(z;\psi^{\pm})$ in \eqref{3BzZ} we obtain
\begin{equation*}\label{St51}
\begin{split}
B(z;\psi_{11})&=\frac{2\alpha}{\alpha\Sigma^{\lambda,\mu}_+(z)+\mu\Sigma^{\alpha,\beta}_+(z)+(\mu-\alpha)z+\alpha\lambda-\beta\mu},\\
B(z;\psi_{12})&=\frac{\beta}{\mu}\left[\frac{\lambda+\mu-z-\Sigma^{\lambda,\mu}_+(z)}{\beta\Sigma^{\lambda,\mu}_+(z)+\lambda\Sigma^{\alpha,\beta}_+(z)+(\beta-\lambda)z+\alpha\lambda-\beta\mu}\right],\\
B(z;\psi_{22})&=\frac{\beta}{\mu}\left[\frac{2\lambda}{\beta\Sigma^{\lambda,\mu}_+(z)+\lambda\Sigma^{\alpha,\beta}_+(z)+(\beta-\lambda)z+\alpha\lambda-\beta\mu}\right].
\end{split}
\end{equation*}
These expressions may be useful to compute the absolutely continuous part of the spectral matrix, but it is difficult to locate the corresponding poles. Therefore, after rationalizing, we obtain
\begin{equation*}
B(z;\psi_{kl})=\frac{p_{kl}(z)+q_{kl}(z)\Sigma^{\lambda,\mu}_+(z)+r_{kl}(z)\Sigma^{\alpha,\beta}_+(z)+s_{kl}(z)\Sigma^{\lambda,\mu}_+(z)\Sigma^{\alpha,\beta}_+(z)}{D(z)},\quad k,l=1,2,
\end{equation*}
where 
\begin{align*}
p_{11}(z)&=((\mu-\alpha)z +\alpha\lambda- \beta\mu)(-z^2+(\alpha + \beta + \lambda + \mu)z+(\lambda - \mu)(\alpha - \beta)),\\
q_{11}(z)&=(\alpha - \mu)z^2+(-\alpha^2 - \alpha\beta - \alpha\lambda + \alpha\mu + 2\beta\mu)z-(\alpha - \beta)(\alpha\lambda - \beta\mu),\\
r_{11}(z)&=(\mu-\alpha)z^2+(2\alpha\lambda + \alpha\mu - \beta\mu - \lambda\mu - \mu^2)z-(\lambda - \mu)(\alpha\lambda - \beta\mu),\\
s_{11}(z)&=(\mu-\alpha)z+\alpha\lambda -\beta\mu,\\
p_{12}(z)&=(\alpha\lambda-2\alpha\beta-2\lambda\mu+ 3\beta\mu)z^2-(\alpha - 3\beta + \lambda - 3\mu)(\alpha\lambda - \beta\mu)z+(\lambda - \mu)(\alpha - \beta)(\alpha\lambda - \beta\mu),\\
q_{12}(z)&=(\alpha\lambda + \beta\mu-2\alpha\beta)z-(\alpha - \beta)(\alpha\lambda - \beta\mu),\\
r_{12}(z)&=(\alpha\lambda + \beta\mu - 2\lambda\mu)z-(\lambda - \mu)(\alpha\lambda - \beta\mu),\\
s_{12}(z)&=\alpha\lambda - \beta\mu,\\
p_{22}(z)&=((\beta-\lambda)z +\alpha\lambda- \beta\mu)(-z^2+(\alpha + \beta + \lambda + \mu)z+(\lambda - \mu)(\alpha - \beta)),\\
q_{22}(z)&=(\beta-\lambda)z^2+(-\alpha\beta + 2\alpha\lambda - \beta^2 + \beta\lambda - \mu)z-(\alpha - \beta)(\alpha\lambda - \beta\mu),\\
r_{22}(z)&=(\lambda-\beta)z^2+(-\alpha\lambda + \beta\lambda + 2\beta\mu - \lambda^2 - \lambda\mu)z-(\lambda - \mu)(\alpha\lambda - \beta\mu),\\
s_{22}(z)&=(\beta-\lambda)z+\alpha\lambda -\beta\mu,\\
D(z)&=4\mu z[(\lambda-\beta)(\alpha-\mu)-(\lambda-\mu+\alpha-\beta)(\alpha\lambda-\beta\mu)z].
\end{align*}
Again, the spectral matrix $\Psi(x)=\Psi_c(x)+\Psi_d(x)$ will have an absolutely continuous part $\Psi_c(x)$ and a discrete part $\Psi_d(x)$. The absolutely continuous part $\Psi_c(x)$ will depend on the position of the closed intervals formed by the zeros of the polynomial inside the square root in \eqref{Sigmab} for $\Sigma^{\lambda,\mu}_{\pm}(z)$ and $\Sigma^{\alpha,\beta}_{\pm}(z)$. Let us call these zeros
$$
\sigma_{\pm}=(\sqrt{\lambda}\pm\sqrt{\mu})^2,\quad \tau_{\pm}=(\sqrt{\alpha}\pm\sqrt{\beta})^2.
$$
We will have 3 different cases, with two sub-cases each:
\begin{enumerate}
\item $[\sigma_-,\sigma_+]\cap[\tau_-,\tau_+]=\emptyset$. We have two situations:
\begin{enumerate}
\item If $\sigma_-<\tau_-$, then the absolutely continuous part $\Psi_c(x)$ is given by
\begin{equation*}
\Psi_c(x)=\begin{cases}
-\D\frac{\Sigma^{\lambda,\mu}_{-}(x)}{\pi D(x)}
\begin{pmatrix}q_{11}(x)+s_{11}(x)\Sigma^{\alpha,\beta}_{+}(x)&q_{12}(x)+s_{12}(x)\Sigma^{\alpha,\beta}_{+}(x)\\q_{12}(x)+s_{12}(x)\Sigma^{\alpha,\beta}_{+}(x)&q_{22}(x)+s_{22}(x)\Sigma^{\alpha,\beta}_{+}(x)\end{pmatrix},\quad x\in[\sigma_-,\sigma_+],\\
\\
-\D\frac{\Sigma^{\alpha,\beta}_{-}(x)}{\pi D(x)}
\begin{pmatrix}r_{11}(x)-s_{11}(x)\Sigma^{\lambda,\mu}_{+}(x)&r_{12}(x)-s_{12}(x)\Sigma^{\lambda,\mu}_{+}(x)\\r_{12}(x)-s_{12}(x)\Sigma^{\lambda,\mu}_{+}(x)&r_{22}(x)-s_{22}(x)\Sigma^{\lambda,\mu}_{+}(x)\end{pmatrix},\quad x\in[\tau_-,\tau_+].
\end{cases}
\end{equation*}
\item If $\tau_-<\sigma_-$, then the absolutely continuous part $\Psi_c(x)$ is given by
\begin{equation*}
\Psi_c(x)=\begin{cases}
-\D\frac{\Sigma^{\alpha,\beta}_{-}(x)}{\pi D(x)}
\begin{pmatrix}r_{11}(x)+s_{11}(x)\Sigma^{\lambda,\mu}_{+}(x)&r_{12}(x)+s_{12}(x)\Sigma^{\lambda,\mu}_{+}(x)\\r_{12}(x)+s_{12}(x)\Sigma^{\lambda,\mu}_{+}(x)&r_{22}(x)+s_{22}(x)\Sigma^{\lambda,\mu}_{+}(x)\end{pmatrix},\quad x\in[\tau_-,\tau_+],\\
\\
-\D\frac{\Sigma^{\lambda,\mu}_{-}(x)}{\pi D(x)}
\begin{pmatrix}q_{11}(x)-s_{11}(x)\Sigma^{\alpha,\beta}_{+}(x)&q_{12}(x)-s_{12}(x)\Sigma^{\alpha,\beta}_{+}(x)\\q_{12}(x)-s_{12}(x)\Sigma^{\alpha,\beta}_{+}(x)&q_{22}(x)-s_{22}(x)\Sigma^{\alpha,\beta}_{+}(x)\end{pmatrix},\quad x\in[\sigma_-,\sigma_+].
\end{cases}
\end{equation*}
\end{enumerate}
\item One interval is strictly contained in the other. We have two situations:
\begin{enumerate}
\item $[\tau_-,\tau_+]\subset[\sigma_-,\sigma_+]$. The absolutely continuous part $\Psi_c(x)$ is given by
\begin{equation*}
\Psi_c(x)=\begin{cases}
-\D\frac{\Sigma^{\lambda,\mu}_{-}(x)}{\pi D(x)}
\begin{pmatrix}q_{11}(x)+s_{11}(x)\Sigma^{\alpha,\beta}_{+}(x)&q_{12}(x)+s_{12}(x)\Sigma^{\alpha,\beta}_{+}(x)\\q_{12}(x)+s_{12}(x)\Sigma^{\alpha,\beta}_{+}(x)&q_{22}(x)+s_{22}(x)\Sigma^{\alpha,\beta}_{+}(x)\end{pmatrix},\quad x\in[\sigma_-,\tau_-],\\
\\
-\D\frac{1}{\pi D(x)}
\begin{pmatrix}q_{11}(x)\Sigma^{\lambda,\mu}_{-}(x)+r_{11}(x)\Sigma^{\alpha,\beta}_{-}(x)&q_{12}(x)\Sigma^{\lambda,\mu}_{-}(x)+r_{12}(x)\Sigma^{\alpha,\beta}_{-}(x)\\q_{12}(x)\Sigma^{\lambda,\mu}_{-}(x)+r_{12}(x)\Sigma^{\alpha,\beta}_{-}(x)&q_{22}(x)\Sigma^{\lambda,\mu}_{-}(x)+r_{22}(x)\Sigma^{\alpha,\beta}_{-}(x)\end{pmatrix},\quad x\in[\tau_-,\tau_+],\\
\\
-\D\frac{\Sigma^{\lambda,\mu}_{-}(x)}{\pi D(x)}
\begin{pmatrix}q_{11}(x)-s_{11}(x)\Sigma^{\alpha,\beta}_{+}(x)&q_{12}(x)-s_{12}(x)\Sigma^{\alpha,\beta}_{+}(x)\\q_{12}(x)-s_{12}(x)\Sigma^{\alpha,\beta}_{+}(x)&q_{22}(x)-s_{22}(x)\Sigma^{\alpha,\beta}_{+}(x)\end{pmatrix},\quad x\in[\tau_+,\sigma_+].
\end{cases}
\end{equation*}

\item $[\sigma_-,\sigma_+]\subset[\tau_-,\tau_+]$. The absolutely continuous part $\Psi_c(x)$ is given by
\begin{equation*}
\Psi_c(x)=\begin{cases}
-\D\frac{\Sigma^{\alpha,\beta}_{-}(x)}{\pi D(x)}
\begin{pmatrix}r_{11}(x)+s_{11}(x)\Sigma^{\lambda,\mu}_{+}(x)&r_{12}(x)+s_{12}(x)\Sigma^{\lambda,\mu}_{+}(x)\\r_{12}(x)+s_{12}(x)\Sigma^{\lambda,\mu}_{+}(x)&r_{22}(x)+s_{22}(x)\Sigma^{\lambda,\mu}_{+}(x)\end{pmatrix},\quad x\in[\tau_-,\sigma_-],\\
\\
-\D\frac{1}{\pi D(x)}
\begin{pmatrix}q_{11}(x)\Sigma^{\lambda,\mu}_{-}(x)+r_{11}(x)\Sigma^{\alpha,\beta}_{-}(x)&q_{12}(x)\Sigma^{\lambda,\mu}_{-}(x)+r_{12}(x)\Sigma^{\alpha,\beta}_{-}(x)\\q_{12}(x)\Sigma^{\lambda,\mu}_{-}(x)+r_{12}(x)\Sigma^{\alpha,\beta}_{-}(x)&q_{22}(x)\Sigma^{\lambda,\mu}_{-}(x)+r_{22}(x)\Sigma^{\alpha,\beta}_{-}(x)\end{pmatrix},\quad x\in[\sigma_-,\sigma_+],\\
\\
-\D\frac{\Sigma^{\alpha,\beta}_{-}(x)}{\pi D(x)}
\begin{pmatrix}r_{11}(x)-s_{11}(x)\Sigma^{\lambda,\mu}_{+}(x)&r_{12}(x)-s_{12}(x)\Sigma^{\lambda,\mu}_{+}(x)\\r_{12}(x)-s_{12}(x)\Sigma^{\lambda,\mu}_{+}(x)&r_{22}(x)-s_{22}(x)\Sigma^{\lambda,\mu}_{+}(x)\end{pmatrix},\quad x\in[\sigma_+,\tau_+].
\end{cases}
\end{equation*}
\end{enumerate}
\item Any other case. We have two situations:
\begin{enumerate}
\item $\sigma_-<\tau_-<\sigma_+<\tau_+$. The absolutely continuous part $\Psi_c(x)$ is given by
\begin{equation*}
\Psi_c(x)=\begin{cases}
-\D\frac{\Sigma^{\lambda,\mu}_{-}(x)}{\pi D(x)}
\begin{pmatrix}q_{11}(x)+s_{11}(x)\Sigma^{\alpha,\beta}_{+}(x)&q_{12}(x)+s_{12}(x)\Sigma^{\alpha,\beta}_{+}(x)\\q_{12}(x)+s_{12}(x)\Sigma^{\alpha,\beta}_{+}(x)&q_{22}(x)+s_{22}(x)\Sigma^{\alpha,\beta}_{+}(x)\end{pmatrix},\quad x\in[\sigma_-,\tau_-],\\
\\
-\D\frac{1}{\pi D(x)}
\begin{pmatrix}q_{11}(x)\Sigma^{\lambda,\mu}_{-}(x)+r_{11}(x)\Sigma^{\alpha,\beta}_{-}(x)&q_{12}(x)\Sigma^{\lambda,\mu}_{-}(x)+r_{12}(x)\Sigma^{\alpha,\beta}_{-}(x)\\q_{12}(x)\Sigma^{\lambda,\mu}_{-}(x)+r_{12}(x)\Sigma^{\alpha,\beta}_{-}(x)&q_{22}(x)\Sigma^{\lambda,\mu}_{-}(x)+r_{22}(x)\Sigma^{\alpha,\beta}_{-}(x)\end{pmatrix},\quad x\in[\tau_-,\sigma_+],\\
\\
-\D\frac{\Sigma^{\alpha,\beta}_{-}(x)}{\pi D(x)}
\begin{pmatrix}r_{11}(x)-s_{11}(x)\Sigma^{\lambda,\mu}_{+}(x)&r_{12}(x)-s_{12}(x)\Sigma^{\lambda,\mu}_{+}(x)\\r_{12}(x)-s_{12}(x)\Sigma^{\lambda,\mu}_{+}(x)&r_{22}(x)-s_{22}(x)\Sigma^{\lambda,\mu}_{+}(x)\end{pmatrix},\quad x\in[\sigma_+,\tau_+].
\end{cases}
\end{equation*}
\item $\tau_-<\sigma_-<\tau_+<\sigma_+$. The absolutely continuous part $\Psi_c(x)$ is given by
\begin{equation*}
\Psi_c(x)=\begin{cases}
-\D\frac{\Sigma^{\alpha,\beta}_{-}(x)}{\pi D(x)}
\begin{pmatrix}r_{11}(x)+s_{11}(x)\Sigma^{\lambda,\mu}_{+}(x)&r_{12}(x)+s_{12}(x)\Sigma^{\lambda,\mu}_{+}(x)\\r_{12}(x)+s_{12}(x)\Sigma^{\lambda,\mu}_{+}(x)&r_{22}(x)+s_{22}(x)\Sigma^{\lambda,\mu}_{+}(x)\end{pmatrix},\quad x\in[\tau_-,\sigma_-],\\
\\
-\D\frac{1}{\pi D(x)}
\begin{pmatrix}q_{11}(x)\Sigma^{\lambda,\mu}_{-}(x)+r_{11}(x)\Sigma^{\alpha,\beta}_{-}(x)&q_{12}(x)\Sigma^{\lambda,\mu}_{-}(x)+r_{12}(x)\Sigma^{\alpha,\beta}_{-}(x)\\q_{12}(x)\Sigma^{\lambda,\mu}_{-}(x)+r_{12}(x)\Sigma^{\alpha,\beta}_{-}(x)&q_{22}(x)\Sigma^{\lambda,\mu}_{-}(x)+r_{22}(x)\Sigma^{\alpha,\beta}_{-}(x)\end{pmatrix},\quad x\in[\sigma_-,\tau_+],\\
\\
-\D\frac{\Sigma^{\lambda,\mu}_{-}(x)}{\pi D(x)}
\begin{pmatrix}q_{11}(x)-s_{11}(x)\Sigma^{\alpha,\beta}_{+}(x)&q_{12}(x)-s_{12}(x)\Sigma^{\alpha,\beta}_{+}(x)\\q_{12}(x)-s_{12}(x)\Sigma^{\alpha,\beta}_{+}(x)&q_{22}(x)-s_{22}(x)\Sigma^{\alpha,\beta}_{+}(x)\end{pmatrix},\quad x\in[\tau_+,\sigma_+].
\end{cases}
\end{equation*}
\end{enumerate}
\end{enumerate}
As for the discrete part $\Psi_d(x)$ we need to study the poles of $B(z;\psi_{kl}), k,l=1,2,$ which are the roots of the second-degree polynomial $D(z)$. One root is 0 and the other the constant
$$
\zeta=\frac{(\alpha-\beta+\lambda-\mu)(\alpha\lambda-\beta\mu)}{(\lambda-\beta)(\alpha-\mu)}.
$$
Defining the constants
$$
C_1=\lambda(\alpha-\mu)^2-\mu(\beta-\lambda)^2,\quad C_2=\beta(\alpha-\mu)^2-\alpha(\beta-\lambda)^2,
$$
we have that the discrete part $\Psi_d(x)$ is given by
\begin{align}
\label{2mm2d}\Psi_d(x)&=\frac{(\beta-\alpha)(\mu-\lambda)}{\mu(\beta-\alpha+\mu-\lambda)}\begin{pmatrix}1&1\\1&1\end{pmatrix}\delta_{0}(x)\mathbf{1}_{\{\mu>\lambda,\beta>\alpha\}}\\
\nonumber&\quad-\frac{C_1C_2}{\mu(\alpha-\mu)(\beta-\lambda)(\alpha-\beta+\lambda-\mu)}\begin{pmatrix}\D\frac{1}{(\beta-\lambda)^2}&\D\frac{1}{(\alpha-\mu)(\beta-\lambda)}\\\D\frac{1}{(\alpha-\mu)(\beta-\lambda)}&\D\frac{1}{(\alpha-\mu)^2}\end{pmatrix}\delta_{\zeta}(x)\mathbf{1}_{\{A_1\cup A_2\cup A_3\}},
\end{align}
where $\mathbf{1}_A$ is the indicator function, $\delta_a(x)$ is the Dirac delta located at $x=a$ and
\begin{align*}
A_1&=\{\alpha>\mu,\;\beta>\lambda,\; C_1<0,\; C_2>0\},\\
A_2&=\{\alpha<\mu,\;\beta>\lambda,\; C_1<0,\; C_2>0\},\\
A_3&=\{\alpha<\mu,\;\beta<\lambda,\; C_1<0,\; C_2>0\}.
\end{align*}
Observe that $A_1\cup A_2\cup A_3$ can also be written as $B\cap(D_1\cup D_2)$ where
$$
B=\left\{\sqrt{\lambda/\mu}<\left|\frac{\beta-\lambda}{\alpha-\mu}\right|<\sqrt{\beta/\alpha}\right\},\quad D_1=\{\beta>\lambda\},\quad D_2=\{\alpha<\mu,\;\beta<\lambda\}.
$$
Finally, the polynomials generated by the three-term recurrence relation \eqref{3TTRRZ} are given by 
\begin{align*}
Q_n^1(x)&=\left(\frac{\mu}{\lambda}\right)^{n/2}U_n\left(\frac{\lambda+\mu-x}{2\sqrt{\lambda\mu}}\right),\quad Q_{-n-1}^1(x)=-\left(\frac{\beta}{\alpha}\right)^{(n+1)/2}U_{n-1}\left(\frac{\alpha+\beta-x}{2\sqrt{\alpha\beta}}\right),\quad n\in\mathbb{N}_0,\\
Q_n^2(x)&=-\left(\frac{\mu}{\lambda}\right)^{(n+1)/2}U_{n-1}\left(\frac{\lambda+\mu-x}{2\sqrt{\lambda\mu}}\right),\quad Q_{-n-1}^2(x)=\left(\frac{\beta}{\alpha}\right)^{n/2}U_{n}\left(\frac{\alpha+\beta-x}{2\sqrt{\alpha\beta}}\right),\quad n\in\mathbb{N}_0,
\end{align*}
where $(U_n)_n$ are the Chebychev polynomials of the second kind. The transition probability functions $P_{ij}(t)$ can then be approximated using the Karlin-McGregor formula \eqref{3FKMcGZ}. It is possible to see that if $\alpha=\mu$ and $\beta=\lambda$ we go back to the Example \ref{exb1} and if $\alpha=\lambda$ and $\beta=\mu$ we go back to the Example \ref{exb2}.

\begin{figure*}[t!]
    \centering
    \begin{subfigure}[b]{0.5\textwidth}
        \centering
        \includegraphics[height=2.5in]{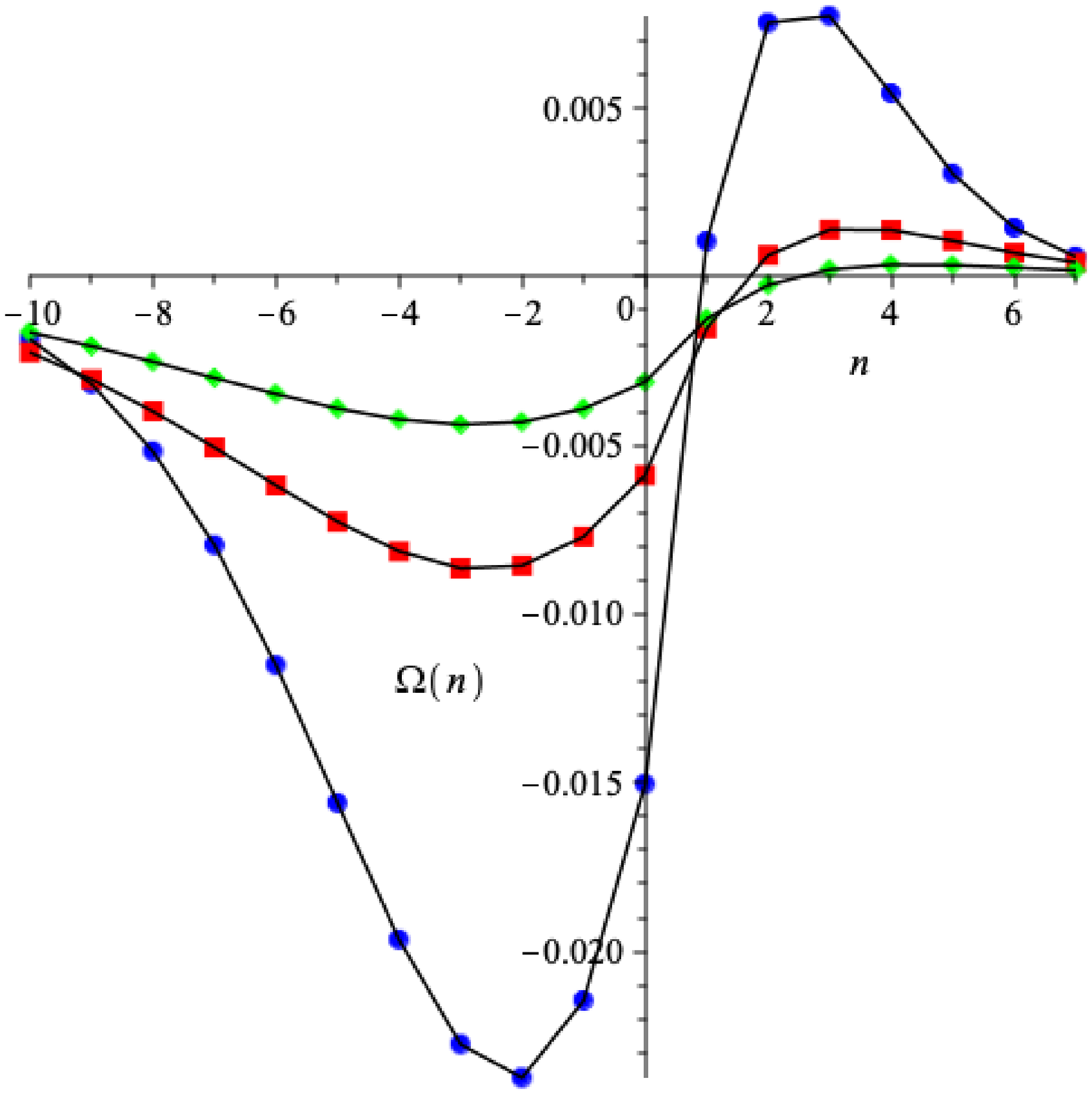}
        \caption{$\lambda=1,\mu=2,\alpha=3,\beta=4$}
    \end{subfigure}%
    ~ 
    \begin{subfigure}[b]{0.5\textwidth}
        \centering
        \includegraphics[height=2.5in]{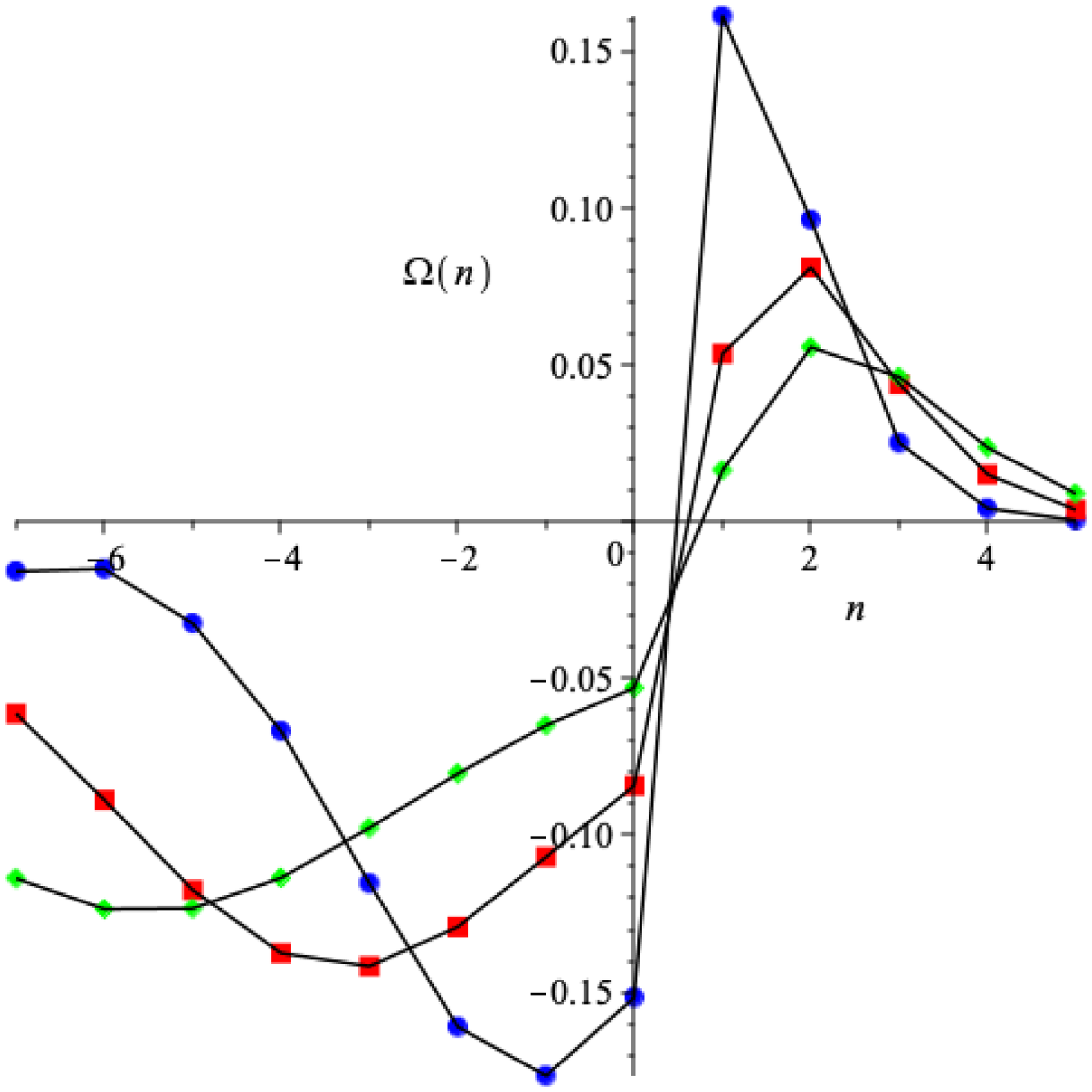}
        \caption{$\lambda=1/2,\mu=1/3,\alpha=13/5,\beta=1/10$}
    \end{subfigure}
    \\
      \centering
    \begin{subfigure}[b]{0.5\textwidth}
        \centering
        \includegraphics[height=2.5in]{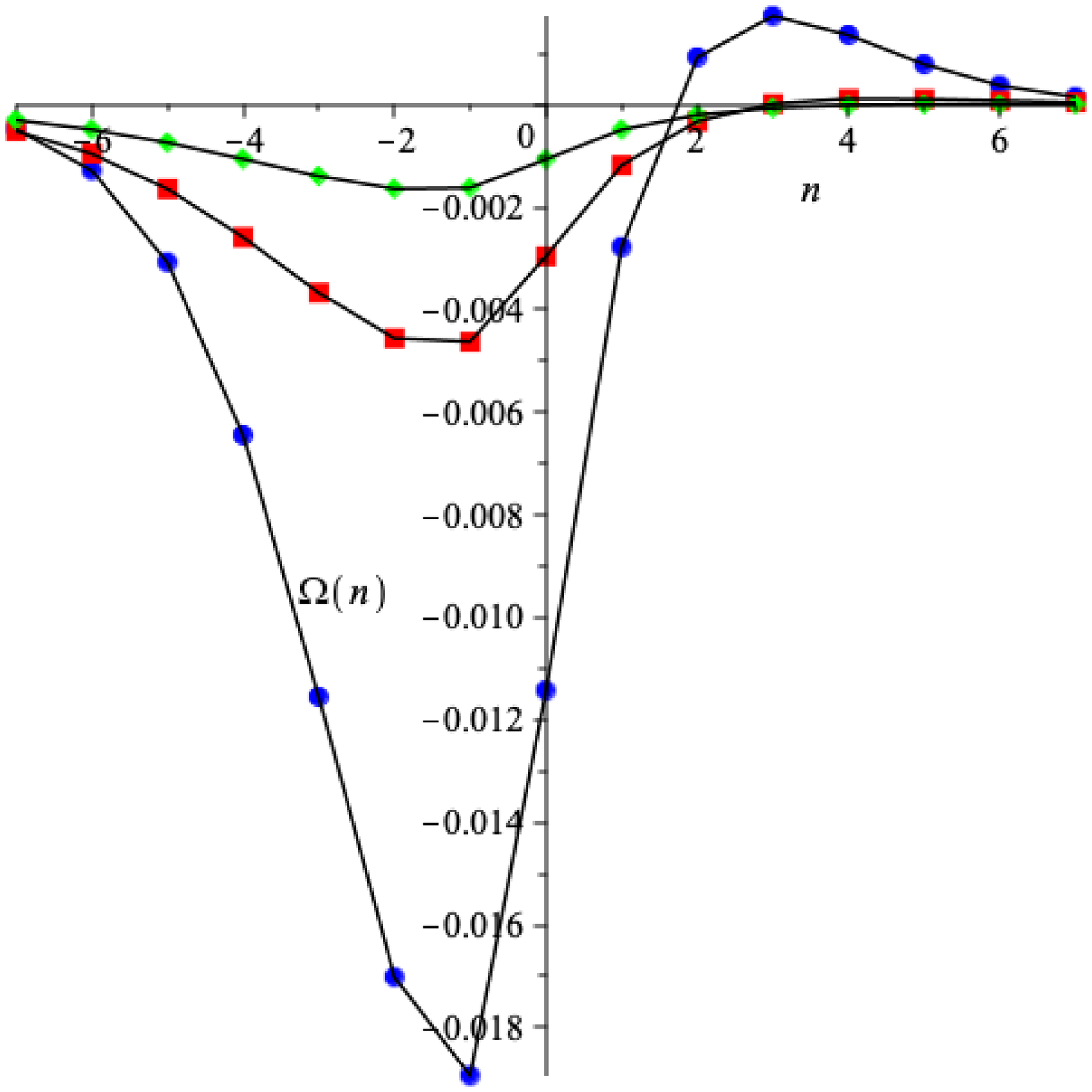}
        \caption{$\lambda=1,\mu=5/2,\alpha=1,\beta=2$}
    \end{subfigure}%
    ~ 
    \begin{subfigure}[b]{0.5\textwidth}
        \centering
        \includegraphics[height=2.5in]{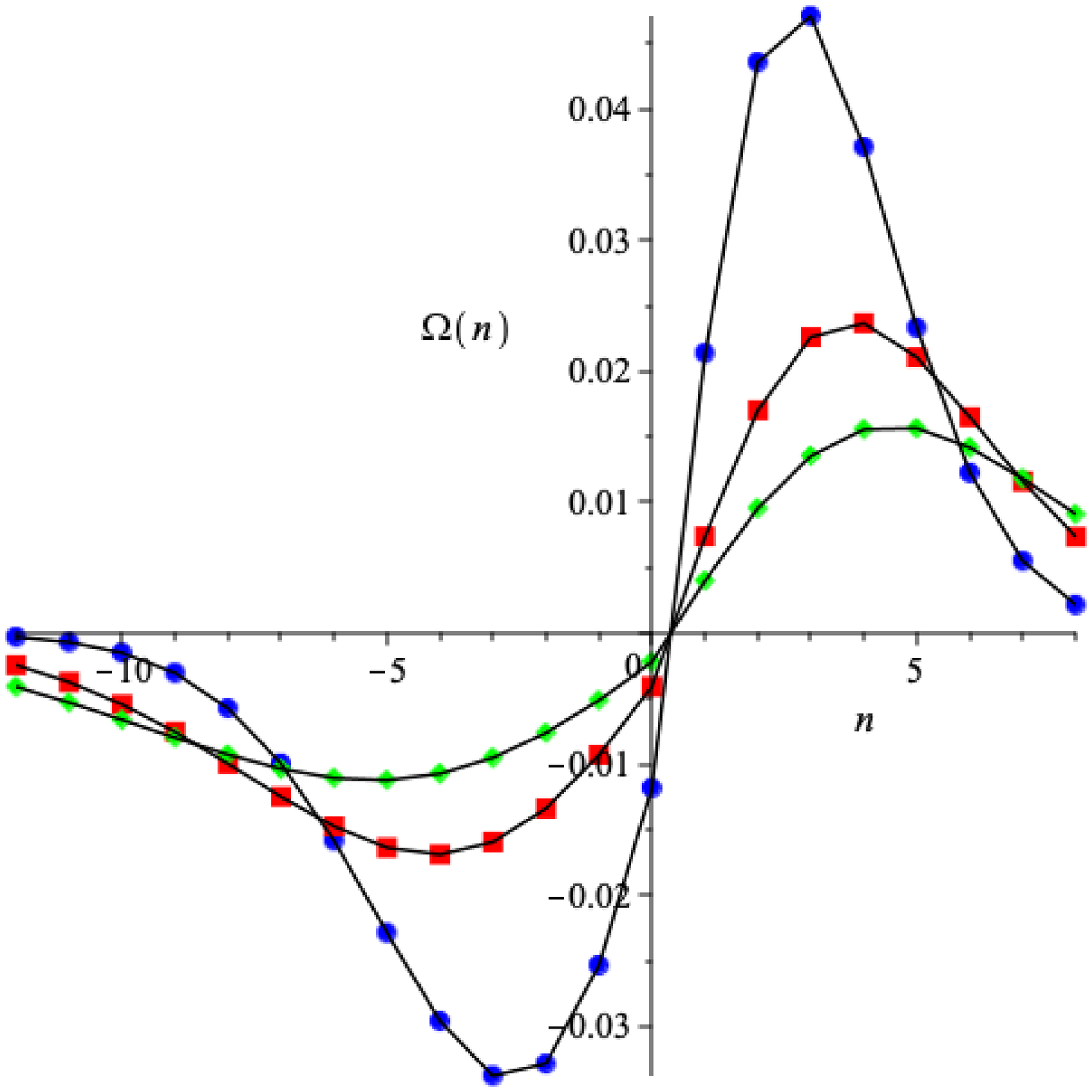}
        \caption{$\lambda=1,\mu=1,\alpha=2,\beta=2$}
    \end{subfigure}
      \caption{The probability current \eqref{probcurZ} (here denoted by $\Omega(n)$) is plotted as function of $n$ for $j = 0$, $t = 3$ (blue circles), $t = 6$ (red squares) and $t =9$ (green diamonds).}  \label{fig7}
\end{figure*}

In all the situations we always have that $\int_0^\infty x^{-1}\psi_{11}(x)dx<\infty$ and $\int_0^\infty x^{-1}\psi_{22}(x)dx<\infty$ for $\lambda>\mu$ \emph{or} $\alpha>\beta$, so the process will be transient. If $\lambda\leq\mu$ \emph{and} $\alpha\leq\beta$ then the integral will diverge and the process will be recurrent. For $\lambda<\mu$ and $\alpha<\beta$ we always have a jump at the point 0, so the process is positive recurrent. If $\lambda=\mu$ and $\alpha=\beta$ then the process is null recurrent. Since the potential coefficients are given here by (see \eqref{3potcoefZ})
$$
\pi_n=\left(\frac{\lambda}{\mu}\right)^{n},\quad n\in\mathbb{N}_0,\quad\pi_{-n}=\frac{\mu}{\beta}\left(\frac{\alpha}{\beta}\right)^{n-1},\quad n\in\mathbb{N},
$$
we have that the invariant distribution $\bm\pi$ for this process is given by
$$
\bm\pi=\frac{(\beta-\alpha)(\mu-\lambda)}{\mu(\beta-\alpha+\mu-\lambda)}\left(\cdots,\frac{\mu\alpha}{\beta^2},\frac{\mu}{\beta},1,\frac{\lambda}{\mu},\frac{\lambda^2}{\mu^2},\cdots\right),\quad \mu>\lambda,\quad\beta>\alpha.
$$
Finally, we can get an approximation of the probability current \eqref{probcurZ}. In Figure \ref{fig7} this probability current is plotted as a function of $n$ starting at $j=0$ for $t=3,6,9$ and for different values of the birth-death rates $\lambda,\mu,\alpha,\beta$. In the first plot we are in the situation $[\sigma_-,\sigma_+]\subset[\tau_-,\tau_+]$ studied in Case (2)(b) and also we have a jump at the point 0 (see \eqref{2mm2d}). In the second case we are in the situation $[\sigma_-,\sigma_+]\cap[\tau_-,\tau_+]=\emptyset$ with $\sigma_-<\tau_-$ studied in Case (1)(a), but there are no discrete jumps. In the third plot we are in the situation $\tau_-<\sigma_-<\tau_+<\sigma_+$ studied in Case (3)(b) and also we have two jumps at the point $0$ and $20/3$ (see \eqref{2mm2d}). Finally, in the fourth case, we are again in the Case (2)(b) but with no jumps.

\end{document}